\pdfoutput=1
\documentclass[preprint,11pt]{elsarticle}
\usepackage[a4paper,margin=2.5cm]{geometry}
\usepackage{amsmath,amssymb,amsthm}
\usepackage{graphicx,booktabs,array}
\usepackage[hidelinks]{hyperref}
\usepackage{microtype}
\usepackage[section]{placeins}
\usepackage{combelow}
\biboptions{sort&compress}
\journal{Applied Numerical Mathematics}
\makeatletter
\def\ps@pprintTitle{\let\@oddhead\@empty\let\@evenhead\@empty
  \def\@oddfoot{\hfil\thepage\hfil}\let\@evenfoot\@oddfoot}
\makeatother

\newtheorem{proposition}{Proposition}[section]
\newtheorem{theorem}[proposition]{Theorem}
\newtheorem{lemma}[proposition]{Lemma}
\newtheorem{corollary}[proposition]{Corollary}
\theoremstyle{remark}
\newtheorem{remark}[proposition]{Remark}

\newcommand{\dd}{\mathrm{d}}


\begin{document}

\begin{frontmatter}

\title{Right shape, wrong time: phase errors of ODE solvers in the flame problem}

\author{Imre Boros\corref{cor1}}
\ead{imre.boros@ictp.acad.ro}
\cortext[cor1]{Corresponding author. ORCID: 0000-0003-2568-4640.}

\affiliation{organization={Tiberiu Popoviciu Institute of Numerical Analysis, Romanian Academy},
            addressline={Str. Fântânele nr. 57, ap. 67--68},
            postcode={400320},
            city={Cluj-Napoca},
            country={Romania}}

\begin{abstract}
The flame model $y'=y^2-y^3$, $y(0)=\delta$, is a standard illustration of stiffness: for $\delta=10^{-4}$ the stiff MATLAB solvers \texttt{ode23s} and \texttt{ode15s} need only $99$ and $140$ steps, against $3040$ for \texttt{ode45}. We show that these cheap solutions have the right shape at the wrong time: they place the front $29$ and $329$ time units early, with maximal errors $0.95$ and $1.00$. Because all solutions are time translates of one master profile, the global error of every numerical solution with values in $[\delta,1)$ is exactly a phase error, $y_k=y_\delta(t_k-D_k)$. For fixed-step one-step methods of order $p$ we prove $D=h^p\int_\delta^{y}g+O(h^{p+1})$ with a remainder bounded uniformly in $\delta$, so that the delay grows at most like $h\ln(1/\delta)$ for $p=1$ and stays bounded for $p\ge2$. For any computed solution the delay is exactly the sum of its defects with respect to the exact flow, converted into time; a heuristic model, confirmed by MATLAB runs, explains why it is created in the first few steps and why tightening the tolerances is an expensive cure. Integrating for $1/y$, which runs like a clock in the slow phase, with an absolute tolerance reduces the error by several orders of magnitude, at a median cost close to that of the standard formulation. The mechanism needs only that the equation is autonomous; its general form is the subject of a companion paper.
\end{abstract}


\end{frontmatter}

%=====================================================================
\section{Introduction}

The initial value problem
\begin{equation}\label{eq:flame}
y'(t)=f(y):=y^2-y^3,\qquad y(0)=\delta,\qquad 0\le t\le 2/\delta,
\end{equation}
is a simple model of the ignition of a ball of flame and has become a standard illustration of stiffness \cite{ShampineGladwellThompson2003,ShampineThompson2007,Moler2004,Moler2014}. For small $\delta>0$ the solution stays close to zero for a long time, undergoes a rapid transition near $t\approx1/\delta$, and then remains close to the equilibrium $y=1$, where the problem is stiff \cite{HairerWanner1996}. Its exact solution is known in terms of the Lambert $W$ function \cite{Corless1996,Moler2014}.

Moler \cite{Moler2014} uses \eqref{eq:flame} to define stiffness: ``a problem is stiff if the solution being sought varies slowly, but there are nearby solutions that vary rapidly''. With $\delta=10^{-4}$ and $\mathrm{RelTol}=10^{-4}$ he reports that \texttt{ode45} takes $3040$ steps and $20179$ function evaluations, while the stiff solvers \texttt{ode23s} and \texttt{ode15s} take $99$ and $140$ steps. Comparisons of this kind, which measure cost only, are the usual way the example is presented \cite{ShampineGladwellThompson2003,Moler2004}. Such cost-only comparisons do not address two questions: can the step counts be predicted, and how accurate are the cheap solutions? Nechita \cite{Nechita2016} proved that the time the exact solution needs to pass between two fixed levels does not depend on $\delta$. \"Ozkaya, Gauger and Nemili \cite{OzkayaGaugerNemili2016} used the same equation to show that explicit Euler with large steps produces oscillatory and chaotic iterates.

This paper answers both questions from one structural observation. The equation is autonomous, so every solution of \eqref{eq:flame} is a time translate of one master profile; equivalently, a solution started at a smaller $\delta$ passes through every larger initial value and coincides with the corresponding solution from then on. Our contributions are:

\begin{enumerate}
\item \emph{An exact error representation} (Sections~\ref{sec:structure} and~\ref{sec:phase}). A closed form of the solution in terms of the Wright $\omega$ function and the exact sensitivity of the front position to the initial value. The global error of every numerical solution with values in $[\delta,1)$ is exactly a phase error, a delay $D_k$ of the front, and for fixed-step one-step methods we prove an asymptotic formula for $D_k$ with a remainder that is uniform in $\delta$.
\item \emph{Cost and accuracy of the MATLAB ODE suite} (Section~\ref{sec:adaptive}). The step counts of explicit solvers follow from their stability boundaries. For adaptive codes the delay is exactly the sum of the local errors converted into time (Proposition~\ref{prop:defects}); a heuristic model built on this identity explains why the cheap stiff solutions misplace the front by an amount that grows like $1/\delta$ or faster, and why tightening the tolerances is an expensive cure.
\item \emph{A remedy that needs no knowledge of the exact solution} (Sections~\ref{sec:sundman} and~\ref{sec:recipe}). Absolute control of the time, either by making time a dependent variable or, more simply, by integrating for $1/y$, which runs like a clock in the slow phase. We give a recipe and demonstrate it with MATLAB, SUNDIALS/CVODE, Julia and SciPy.
\item \emph{A second example} (Section~\ref{sec:q3}). The family $y'=y^q(1-y)$, with a complete study of $q=3$.
\end{enumerate}
Sections~\ref{sec:structure} and~\ref{sec:phase} contain proofs. Everything else is either an exact identity (Lemma~\ref{lem:phase}, Proposition~\ref{prop:defects}) evaluated on computed solutions, or a heuristic model that we state as such and test against the computations. The main solver study uses MATLAB R2023b; results obtained with SciPy, SUNDIALS/CVODE and Julia are identified as such. The scripts and data that reproduce every figure and table, including Tables~S1--S6 of the supplementary material, are openly available \cite{BorosCode2026}.

\paragraph{Relation to prior work} The time-shift character of the error is classical for periodic orbits, where the error grows along the orbit while it stays small across it \cite{CanoSanzSerna1997}; for a scalar autonomous problem the decomposition is exact (Lemma~\ref{lem:phase}). The exponents $\alpha=p/(p+1)$ and $\alpha=1$ of our delay law are those of the analysis of tolerance proportionality for error control without and with local extrapolation \cite{Higham1993,SoderlindWang2006}; what is specific here is the factor $\delta^{1-q}$, which comes from the sensitivity of the front position and makes any fixed tolerance useless as $\delta\to0$. Goal-oriented time adaptivity \cite{MeisrimelBirken2020} chooses the steps from their contribution to the error in a quantity of interest given by a time integral of the solution, and Chaudhry, Estep, Stevens and Tavener \cite{ChaudhryEstepStevensTavener2021} estimate the error in the first time at which the solution crosses a threshold, a regular crossing like our front, by adjoint-based a posteriori analysis. Event location \cite{ShampineGladwellThompson2003} finds the time at which the computed solution crosses a level, but it cannot remove a delay that is already in the computed solution (Section~\ref{sec:adaptive}). Sundman transformations have been used to compute blow-up times \cite{StuartFloater1990}, and Hoel and Meo \cite{HoelMeo2026} choose step sizes a priori from the sensitivity of the blow-up time. Higham, Humphries and Wain \cite{HighamHumphriesWain2000} observed that standard error control performs poorly near a fixed point, where orbits are bounded in space but unbounded in time, and proposed to bound the local error by a fraction of the arc length of the step; with a related constraint Humphries and Vigneswaran \cite{HumphriesVigneswaran2004} removed a phase shift in the passage near a hyperbolic saddle. Their analysis is linear, for hyperbolic fixed points, and does not consider the error of an event time; the equilibrium $y=0$ of $y'=y^q(1-y)$, $q\ge2$, is degenerate, and there such a constraint gives a relative accuracy in time and no more, so that for the flame problem the delay in units of the front width still grows like $1/\delta$ \cite{BorosC2026}. Our approach differs from these in method rather than in target: for a scalar autonomous problem the delay is available exactly through the time function $T_\delta$, the fixed-step analysis gives a remainder that is uniform in the degeneracy parameter $\delta$, and the remedy is a change of unknown that lets the error test of an unmodified standard solver measure time. In a companion preprint \cite{BorosB2026} the master profile is computed by spectral collocation with the front pinned at a prescribed time; there the same sensitivity of the front position appears as an ill-conditioned translation mode of the initial value formulation. A second companion paper \cite{BorosC2026} shows that the mechanism needs only that the equation is autonomous and treats systems and other problems (epidemic models, ignition, neurons near threshold, a heteroclinic cycle): it defines a condition number of the event time in the error norm of the solver, extends the delay law and gives a detection algorithm. The present paper is the detailed study of the flame family, for which every statement is exact. The notation is summarized in Table~\ref{tab:notation}.

\begin{table}[!tbp]
\centering\small
\caption{Notation.}
\label{tab:notation}
\renewcommand{\arraystretch}{1.15}
\begin{tabular}{>{\raggedright\arraybackslash}p{0.24\textwidth}p{0.68\textwidth}}
\toprule
$y_\delta(t)$, $f$, $f_q$ & exact solution with $y(0)=\delta$; $f(y)=y^2-y^3$; $f_q(y)=y^q(1-y)$, $f=f_2$\\
$T_\delta(y)$, $c(\delta)$ & time the exact solution needs to reach $y$; position of the front, $y_\delta(c(\delta))=1/2$\\
$(t_k,y_k)$, $h$, $h_k$ & computed points (accepted steps); fixed and variable step size\\
$D_k=t_k-T_\delta(y_k)$ & phase delay; $D_k>0$: the computed solution lags (front late), $D_k<0$: front early\\
$D_{\mathrm{front}}$, $|D|$ & delay at the first computed point with $y_k\ge1/2$ (front delay)\\
$e_k$, maximal error & $e_k=y_k-y_\delta(t_k)$; maximal error $=\max_k|e_k|$ over the computed points (maximal nodal error)\\
$p$, $g$, $\kappa$ & order of the method; delay density of Theorem~\ref{prop:delay}; $\kappa=\lim_{y\to0}y\,g(y)$\\
$\mathrm{RelTol}$, $\mathrm{AbsTol}$, tol & tolerances of the error test $|e_i|\le\max(\mathrm{RelTol}\,|u_i|,\mathrm{AbsTol}_i)$; tol: common parameter of a tolerance strategy\\
benchmark tolerance & $\mathrm{RelTol}=10^{-4}$, $\mathrm{AbsTol}=10^{-6}$, the setting of \cite{Moler2014}\\
$\alpha$ & exponent of the tolerance in the delay law \eqref{eq:adaptive}\\
$L_{\mathrm{stiff}}$, $\beta$ & length of the stiff phase; real stability boundary of an explicit method\\
steps, $f$-evals (nfev) & accepted steps; number of function evaluations\\
$u$ & unit roundoff of double precision, $2^{-53}\approx1.1\cdot10^{-16}$\\
$r_{\min}$ & smallest RelTol MATLAB accepts, $100\,\mathrm{eps}=100\cdot2^{-52}\approx2.2\cdot10^{-14}$\\
\bottomrule
\end{tabular}
\end{table}

%=====================================================================
\section{Translation structure of the flame problem}\label{sec:structure}

Let $0<\delta<1$ and $\lambda:=1/\delta-1$. Since $f>0$ on $(0,1)$, the time needed by the solution of \eqref{eq:flame} to reach the level $y\in[\delta,1)$ is
\begin{equation}\label{eq:T}
T_\delta(y)=\int_\delta^y\frac{\dd z}{z^2(1-z)}=\Bigl(\frac1\delta-\frac1y\Bigr)+\ln\frac{y(1-\delta)}{\delta(1-y)} .
\end{equation}
Write $w:=1/y-1$. Then \eqref{eq:T} reads $t=\lambda+\ln\lambda-(w+\ln w)$. The Wright omega function $\omega$ is defined for real arguments by $\omega(z)+\ln\omega(z)=z$ \cite{CorlessJeffrey2002}; it is related to the Lambert function by $\omega(z)=W_0(e^z)$ for real $z$. Hence
\begin{equation}\label{eq:exact}
y_\delta(t)=\frac{1}{1+\omega(\lambda+\ln\lambda-t)} ,
\end{equation}
which is the formula $y=1/(W(\lambda e^{\lambda-t})+1)$ of \cite{Moler2014} written without the factor $e^{\lambda}$. This matters in practice: $\lambda e^{\lambda-t}$ overflows in double precision as soon as $\lambda+\ln\lambda>\ln(1.8\cdot10^{308})\approx709.8$, i.e.\ for $\lambda>703$ or $\delta<1.4\cdot10^{-3}$, whereas \eqref{eq:exact} can be evaluated for any $\delta$ (we use \texttt{scipy.special.wrightomega} \cite{SciPy2020}).

\begin{proposition}[Translation]\label{prop:translation}
Let $Y(s):=1/(1+\omega(1-s))$, so that $Y(0)=1/2$, and $c(\delta):=\lambda+\ln\lambda-1$. Then
\[
y_\delta(t)=Y\bigl(t-c(\delta)\bigr)\qquad\text{for all }t,\ \delta\in(0,1).
\]
In particular, for $0<\delta_1<\delta_2<1$ we have $y_{\delta_1}(t+\tau)=y_{\delta_2}(t)$ with $\tau=c(\delta_1)-c(\delta_2)=T_{\delta_1}(\delta_2)$.
\end{proposition}

\begin{proof}
Immediate from \eqref{eq:exact}; the second statement is $y_{\delta_1}(\tau)=\delta_2$ together with uniqueness for the autonomous equation.
\end{proof}

\begin{figure}[t]
\centering
\includegraphics[width=\textwidth]{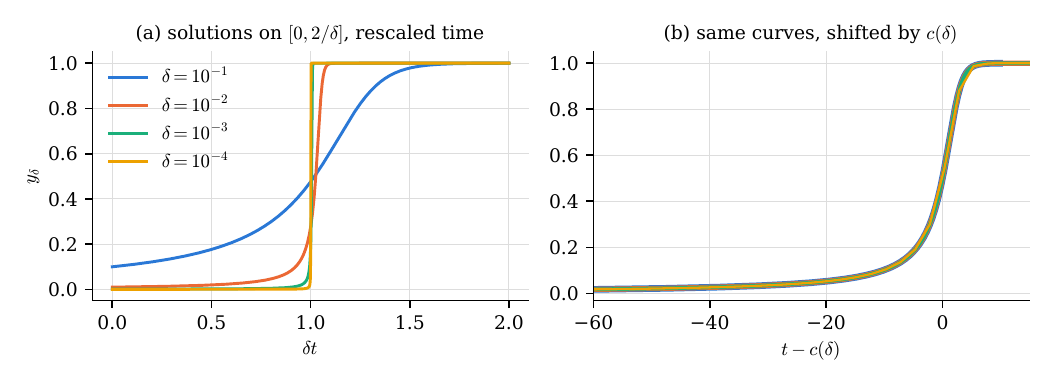}
\caption{(a) Solutions of \eqref{eq:flame} for four values of $\delta$, plotted against $\delta t$. (b) The same solutions plotted against $t-c(\delta)$ collapse onto the master profile $Y$.}
\label{fig:translation}
\end{figure}

The parameter $\delta$ thus affects only \emph{where} the master profile is placed, not its shape (Fig.~\ref{fig:translation}). Expanding $c$,
\begin{equation}\label{eq:c}
c(\delta)=\frac1\delta+\ln\frac1\delta-2+O(\delta),
\end{equation}
so the transition is centered at $1/\delta+\ln(1/\delta)$, not at $1/\delta$; the logarithmic term is the shift, not a width.

\begin{corollary}[Transition times]\label{cor:width}
For $0<\alpha<\beta<1$ the solution passes from $\alpha$ to $\beta$ in time
\[
T_\delta(\beta)-T_\delta(\alpha)=\Bigl(\frac1\alpha-\frac1\beta\Bigr)+\ln\frac{\beta(1-\alpha)}{\alpha(1-\beta)},
\]
independently of $\delta$. For example, the passage from $0.1$ to $0.9$ takes $13.28$ time units for every $\delta<0.1$.
\end{corollary}

This was proved by Nechita \cite[formula (6)]{Nechita2016}. Proposition~\ref{prop:translation} shows that the $\delta$-independence is not specific to two levels: the whole profile is $\delta$-independent.

\begin{proposition}[Sensitivity of the transition time]\label{prop:sens}
If the initial value is perturbed to $\delta(1+\rho)$, $-1<\rho<1/\delta-1$, the solution is translated in time by
\[
T_\delta\bigl(\delta(1+\rho)\bigr)=\frac{\rho}{\delta(1+\rho)}+\ln\frac{(1+\rho)(1-\delta)}{1-\delta(1+\rho)}=\frac\rho\delta+\rho+O(\delta\rho)+O(\rho^2/\delta).
\]
\end{proposition}

\begin{proof}
By Proposition~\ref{prop:translation} the perturbed solution is $y_\delta(t+T_\delta(\delta(1+\rho)))$; insert $y=\delta(1+\rho)$ into \eqref{eq:T}.
\end{proof}

A relative error $\rho$ committed while $y\approx\delta$ therefore moves the whole transition by $\rho/\delta$. This is the precise content of the ``mild instability'' of \eqref{eq:flame} on $[0,1/\delta]$ \cite{ShampineGladwellThompson2003}: in the initial part of the slow phase the linearization $J=f'(y)=2y-3y^2$ is $O(\delta)$, but the slow phase lasts a time $O(1/\delta)$, and the net effect is a pure time shift.

\paragraph{The slow phase as a regularized blow-up} For $y\ll1$ the equation is close to $y'=y^2$, whose solution $\delta/(1-\delta t)$ blows up at $t=1/\delta$. The term $-y^3$ stops the growth at $y=1$, at the price of the shift $\ln(1/\delta)$ in \eqref{eq:c}. In the variables $z=y/\delta$, $s=\delta t$ the problem becomes $\dd z/\dd s=z^2(1-\delta z)$, $z(0)=1$, and \eqref{eq:T} gives
\[
\frac1z=1-s+\delta\ln\frac{z(1-\delta)}{1-\delta z},
\]
so that $z\to1/(1-s)$ as $\delta\to0$, uniformly on $[0,s_0]$ for every $s_0<1$. The slow phase is thus a regularized blow-up of $y'=y^2$, and methods developed for computing blow-up apply to it (Section~\ref{sec:sundman}).

%=====================================================================
\section{Phase errors of one-step methods}\label{sec:phase}

\subsection{Global error equals phase error}

\begin{lemma}\label{lem:phase}
Let $(t_k,y_k)$ be any numerical approximation of \eqref{eq:flame} with $\delta\le y_k<1$, and define the \emph{phase delay}
\[
D_k:=t_k-T_\delta(y_k).
\]
Then $y_k=y_\delta(t_k-D_k)$ exactly, and the global error satisfies
\[
e_k:=y_k-y_\delta(t_k)=-D_k\,f\bigl(y_\delta(\theta_k)\bigr)\quad\text{for some }\theta_k\text{ between }t_k-D_k\text{ and }t_k,
\qquad |e_k|\le\tfrac{4}{27}|D_k|.
\]
\end{lemma}

\begin{proof}
$y_\delta(T_\delta(y_k))=y_k$ by definition of $T_\delta$; then apply the mean value theorem and $\max_{[0,1]}f=f(2/3)=4/27$.
\end{proof}

\begin{figure}[!tbp]
\centering
\includegraphics[width=\textwidth]{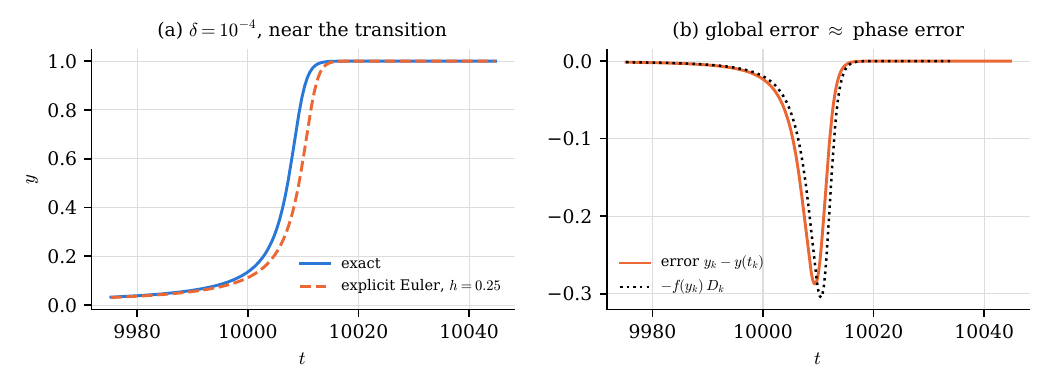}
\caption{Explicit Euler, $h=0.25$, $\delta=10^{-4}$. (a) Exact and numerical solutions near the front. (b) Global error and its first-order phase representation $-f(y_k)D_k$ (Lemma~\ref{lem:phase}).}
\label{fig:errphase}
\end{figure}

For a scalar autonomous problem the numerical solution is therefore the exact solution sampled at wrong times. All information about the error is in the scalar sequence $D_k$, which can be computed exactly from \eqref{eq:T}. Near the front $e_k\approx-f(y_k)D_k$ (Fig.~\ref{fig:errphase}); once the solution has passed the front, $f(y)\to0$ and the error disappears although $D_k$ remains. The phenomenon that the error is dominated by a time shift is well known for periodic orbits \cite{CanoSanzSerna1997}; here it is exact.

In terms of the sign, $D_k>0$ means that the numerical solution lags behind the exact one, and $D_k<0$ that it runs ahead; a front computed `early', as by \texttt{ode15s} in Section~\ref{sec:adaptive}, has $D<0$. The lemma requires only $\delta\le y_k<1$. It does not assume fixed steps, a particular method or monotone iterates; monotonicity is needed only to read $D_k$ as the delay of a front that is crossed once. The identity holds point by point; that the computed solution is close to a single translate of the exact one is an additional, empirical observation: for the MATLAB solvers the delay changes by at most $0.1\%$ across the front (Section~\ref{sec:adaptive}). For explicit Euler, $1-y_{k+1}=(1-y_k)(1-hy_k^2)$ with $0<1-hy_k^2<1$, so for every $h\le1$ the iterates increase monotonically and stay in $[\delta,1)$; other methods have similar method-dependent step-size bounds. For larger steps the iterates can overshoot $y=1$; this is the regime of the oscillatory and chaotic iterates studied in \cite{OzkayaGaugerNemili2016}, where the error is no longer a phase error.

\subsection{A uniform asymptotic formula for the phase delay}

Section~\ref{sec:q3} treats a second example, so we work with the family
\begin{equation}\label{eq:fq}
y'=f_q(y):=y^q(1-y),\qquad y(0)=\delta,\qquad q\in\{2,3,\dots\},
\end{equation}
which contains \eqref{eq:flame} for $q=2$. Partial fractions give $T_\delta(y)=\Theta_q(y)-\Theta_q(\delta)$ with
\begin{equation}\label{eq:Theta}
\Theta_q(y)=\ln\frac{y}{1-y}-\sum_{j=1}^{q-1}\frac{1}{j\,y^{j}} ,
\end{equation}
and Lemma~\ref{lem:phase} holds for \eqref{eq:fq} with $\tfrac4{27}$ replaced by $m_q:=\max_{[0,1]}f_q=q^q/(q+1)^{q+1}$.

Near $y=0$ the equation is invariant under $y\mapsto\lambda y$, $t\mapsto\lambda^{1-q}t$, and the natural step variable is $z:=h\,y^{q-1}$. The following lemma makes this precise for the exact flow $\mathcal{F}_h$.

\begin{lemma}[Scaled exact flow]\label{lem:flow}
Let $W(\sigma,y)$ denote the value at time $\sigma$ of the solution of $w'=w^q(1-yw)$, $w(0)=1$. Then $\mathcal{F}_h(y)=y\,W(hy^{q-1},y)$ for $0<y<1$, $h\ge0$. There is $r_0>0$ such that $W$ is analytic on a complex neighbourhood of $\{|\sigma|\le r_0\}\times[0,1]$, and $W(0,y)=1$.
\end{lemma}

\begin{proof}
If $Y'=Y^q(1-Y)$, $Y(0)=y$, then $w(\sigma):=Y(\sigma y^{1-q})/y$ satisfies $\dd w/\dd\sigma=y^{-q}Y^q(1-Y)=w^q(1-yw)$ and $w(0)=1$. The right-hand side is a polynomial in $(w,y)$, so the solution depends analytically on $(\sigma,y)$ near $\sigma=0$, uniformly for $y$ in a compact complex neighbourhood of $[0,1]$ (Cauchy's majorant method; see e.g.\ \cite{HairerNorsettWanner1993}).
\end{proof}

\noindent\textbf{Assumption (S).} The method has the form $\Phi_h(y)=y\,\psi(hy^{q-1},y)$ with $\psi$ analytic on a complex neighbourhood of $\{|z|\le r\}\times[0,1]$ for some $r>0$, and it has order $p\ge1$ in the sense that
\[
\psi(z,y)-W(z,y)=z^{p+1}E(z,y),
\]
with $E$ analytic on the same set.

\begin{lemma}\label{lem:S}
Assumption (S) holds for every explicit Runge--Kutta method and every implicit Runge--Kutta method of classical order $p$.
\end{lemma}

\begin{proof}
For a Runge--Kutta method write the stages as $Y_i=yW_i$. The stage equations become $W_i=1+z\sum_j a_{ij}W_j^q(1-yW_j)$, and $\Phi_h(y)=y\bigl(1+z\sum_ib_iW_i^q(1-yW_i)\bigr)$. For explicit methods the $W_i$ are polynomials in $(z,y)$. For implicit methods the implicit function theorem at $z=0$, $W_i=1$, gives analytic $W_i$ for $|z|$ small, uniformly on a neighbourhood of $[0,1]$. In both cases $\psi$ is analytic. Expanding in powers of $h$ at fixed $y\in(0,1)$, $\Phi_h(y)-\mathcal{F}_h(y)=\sum_{j\ge0}h^j y^{1+j(q-1)}\bigl[\psi_j(y)-W_j(y)\bigr]$, where $\psi_j,W_j$ are the Taylor coefficients in $z$. Classical order $p$ means that this vanishes to order $h^{p+1}$, hence $\psi_j=W_j$ for $j\le p$ on $(0,1)$ and, by analyticity, identically.
\end{proof}

Under (S) the leading local error is $\Phi_h(y)-\mathcal{F}_h(y)=h^{p+1}C(y)+O(h^{p+2})$ with $C(y)=y^{1+(p+1)(q-1)}E(0,y)$; in particular $C$ vanishes to high order at $y=0$. We define
\[
c(y):=\frac{C(y)}{f_q(y)}=\frac{y^{p(q-1)}E(0,y)}{1-y},\qquad g(y):=-\frac{C(y)}{f_q(y)^2}=-\frac{y^{p(q-1)-q}E(0,y)}{(1-y)^2}.
\]

\begin{theorem}[Uniform phase-delay asymptotics]\label{prop:delay}
Let $q\ge2$, let the method satisfy (S) with order $p$, and let $\eta\in(0,1)$. There are $h_0>0$ and $K>0$, depending on $q$, $\eta$ and the method but \emph{not on $\delta$}, such that for all $h\in(0,h_0]$, all $\delta\in(0,\eta)$ and all $n\ge1$ with $y_0,\dots,y_{n-1}\le\eta$ (that is, up to and including the first step that leaves $[\delta,\eta]$):
\begin{equation}\label{eq:delay}
\Bigl|\,D_n-h^p\int_\delta^{y_n}g(y)\,\dd y\,\Bigr|\le K\,h^{p+1}.
\end{equation}
Moreover the iterates are strictly increasing.
\end{theorem}

\begin{proof}
Throughout, $K$ denotes constants independent of $\delta$, $h$ and $n$. Put $\eta'=(1+\eta)/2$, $\eta''=(1+\eta')/2$ and $z=hy^{q-1}\in(0,h]$.

\emph{Step 1 (one step).} Since $\psi(0,y)=W(0,y)=1$ and both are analytic, $|\psi(z,y)-1|+|W(z,y)-1|\le Kz$ for $y\in[0,1]$, $0\le z\le r$. Hence, for $y\in(0,\eta']$ and $h\le h_0$ small, $\Phi_h(y)$ and $\mathcal{F}_h(y)$ lie in $[y(1-Kz),y(1+Kz)]\subset(0,\eta'']$, and in particular $\Phi_h(y)\le2y$. Define the one-step time defect
\[
\tau(y):=T_\delta(\Phi_h(y))-T_\delta(y)-h=\Theta_q(\Phi_h(y))-\Theta_q(\mathcal{F}_h(y))=\int_{\mathcal{F}_h(y)}^{\Phi_h(y)}\frac{\dd x}{f_q(x)},
\]
where we used $\Theta_q(\mathcal{F}_h(y))=\Theta_q(y)+h$. Note that $\tau$ does not depend on $\delta$. For $x$ between $\mathcal{F}_h(y)$ and $\Phi_h(y)$, $f_q(x)/f_q(y)=(x/y)^q(1-x)/(1-y)=1+O(z)$ because $1-y\ge1-\eta'$. With $\Phi_h(y)-\mathcal{F}_h(y)=y z^{p+1}\bigl(E(0,y)+O(z)\bigr)$ we obtain
\[
\tau(y)=\frac{yz^{p+1}}{f_q(y)}\bigl(E(0,y)+O(z)\bigr)=h^{p+1}c(y)+\rho(y),\qquad |\rho(y)|\le Kh^{p+2}y^{(p+1)(q-1)},
\]
since $yz^{p+1}/f_q(y)=h^{p+1}y^{p(q-1)}/(1-y)$ and $yz^{p+2}/f_q(y)=h^{p+2}y^{(p+1)(q-1)}/(1-y)$. In particular $|\tau(y)|\le Kh^{p+1}y^{p(q-1)}\le h/2$ for $h\le h_0$.

\emph{Step 2 (two clocks).} Let $T_k:=T_\delta(y_k)$. By Step 1, $T_{k+1}-T_k=h+\tau(y_k)\in[h/2,3h/2]$ for $k<n$, so $T_k$, and hence $y_k=y_\delta(T_k)$, is strictly increasing, and
\[
D_n=nh-T_n=-\sum_{k<n}\tau(y_k)=-h^p\sum_{k<n}h\,c(y_k)-\sum_{k<n}\rho(y_k).
\]

\emph{Step 3 (Riemann sum).} Since $c(y)=y^{p(q-1)}e(y)$ with $e$ analytic near $[0,\eta'']$, $|c'(\xi)|\le K\xi^{p(q-1)-1}$, and $p(q-1)-1\ge0$. For $s\in[T_k,T_{k+1}]$ we have $y_\delta(s)\in[y_k,y_{k+1}]\subset[y_k,2y_k]$ and $|y_\delta(s)-y_k|\le\frac32h\max_{[y_k,2y_k]}f_q\le Khy_k^q$. Therefore $|c(y_\delta(s))-c(y_k)|\le Khy_k^{p(q-1)-1+q}=Khy_k^{(p+1)(q-1)}$ and
\[
\Bigl|\int_{T_k}^{T_{k+1}}c(y_\delta(s))\,\dd s-(h+\tau(y_k))\,c(y_k)\Bigr|\le Kh^2y_k^{(p+1)(q-1)}.
\]
Moreover $|\tau(y_k)c(y_k)|\le Kh^{p+1}y_k^{2p(q-1)}\le Kh^2y_k^{(p+1)(q-1)}$ because $p\ge1$ and $y_k\le1$. Hence
\begin{equation}\label{eq:riemann}
\Bigl|h\,c(y_k)-\int_{T_k}^{T_{k+1}}c(y_\delta(s))\,\dd s\Bigr|\le Kh^2y_k^{(p+1)(q-1)}.
\end{equation}

\emph{Step 4 (a sum bounded independently of $\delta$).} Since $T_{k+1}-T_k\ge h/2$ and $y_\delta\ge y_k$ on $[T_k,T_{k+1}]$,
\[
h\sum_{k<n}y_k^{(p+1)(q-1)}\le2\int_0^{T_n}y_\delta(s)^{(p+1)(q-1)}\dd s=2\int_\delta^{y_n}\frac{y^{p(q-1)-1}}{1-y}\,\dd y\le\frac{2}{p(q-1)(1-\eta'')}=:S .
\]
The exponent $p(q-1)-1$ is nonnegative, so the integral is bounded independently of $\delta$. This is the crucial point: the long slow phase near $y=\delta$ contributes only a bounded amount.

\emph{Step 5 (conclusion).} Summing \eqref{eq:riemann} over $k<n$ and using
\[
\int_0^{T_n}c(y_\delta(s))\,\dd s=\int_\delta^{y_n}\frac{c(y)}{f_q(y)}\,\dd y=-\int_\delta^{y_n}g(y)\,\dd y
\]
gives $\bigl|\sum_{k<n}hc(y_k)+\int_\delta^{y_n}g\bigr|\le KhS$. Together with $\sum_{k<n}|\rho(y_k)|\le Kh^{p+1}S$ and Step 2 this yields \eqref{eq:delay}.
\end{proof}

\begin{corollary}\label{cor:delay}
Under the assumptions of Theorem~\ref{prop:delay}, uniformly in $\delta\in(0,\eta)$:
\begin{enumerate}
\item if $p=1$: $D_n=\kappa\,h\ln(y_n/\delta)+O(h)$ with $\kappa=-E(0,0)=\lim_{y\to0}y\,g(y)$;
\item if $p\ge2$: $D_n=h^p\int_0^{y_n}g(y)\,\dd y+O\bigl(h^p\delta^{\,p(q-1)-q+1}+h^{p+1}\bigr)$, a delay independent of $\delta$ to leading order;
\item the global error satisfies $|e_n|\le m_q|D_n|\le K\,h^p\ln(1/\delta)$ for $p=1$ and $\le K\,h^p$ for $p\ge2$, as long as $y_{n-1}\le\eta$.
\end{enumerate}
\end{corollary}

\begin{proof}
$g(y)=-y^{p(q-1)-q}E(0,y)/(1-y)^2$. For $p=1$ the exponent is $-1$; for $p\ge2$ it is $\ge q-2\ge0$. Item~3 follows from Lemma~\ref{lem:phase}.
\end{proof}

\begin{remark}[Scope of Theorem~\ref{prop:delay}]\label{rem:scope}
The theorem concerns fixed step sizes and the part of the trajectory with $y\le\eta$; the constants $h_0$ and $K$ depend on $\eta$ and deteriorate as $\eta\to1$. It says nothing about the stiff phase $y\approx1$, where $f_q(y)\to0$ and the error disappears although $D_k$ remains (Lemma~\ref{lem:phase}). For implicit methods, Assumption~(S) includes the solvability of the stage equations, which holds for $h\le h_0$ by the implicit function theorem (proof of Lemma~\ref{lem:S}); $h_0$ does not depend on $\delta$ because the stage equations are posed in the scaled variable $z=hy^{q-1}\le h$. Adaptive step sizes are not covered: for them we use the exact identity of Proposition~\ref{prop:defects} below, and the delay law derived from it is a model, not a theorem.
\end{remark}

For comparison, the classical global error bound for a method of order $p$ \cite{HairerNorsettWanner1993} contains the factor $\exp\bigl(\int_0^{t_n}f_q'(y_\delta)\,\dd t\bigr)=f_q(y_n)/f_q(\delta)\sim\delta^{-q}$ coming from the instability of the slow phase. Corollary~\ref{cor:delay} shows that this factor is an artefact of the norm-wise analysis: measured as a phase, the error does not grow at all with the length of the slow phase for $p\ge2$, and only logarithmically for $p=1$. The leading term of \eqref{eq:delay} is consistent with the classical asymptotic expansion of the global error \cite[Sect.~II.8]{HairerNorsettWanner1993}: for fixed $\delta$, $e_n=h^pe_p(t_n)+O(h^{p+1})$, where $e_p$ solves the variational equation $e_p'=f_q'(y_\delta)e_p+C(y_\delta)$, $e_p(0)=0$. For a scalar autonomous problem this equation can be integrated, $e_p(t)=f_q(y_\delta(t))\int_0^tC(y_\delta(s))/f_q(y_\delta(s))\,\dd s$, and $-e_p/f_q$ is the leading term of \eqref{eq:delay}. What Theorem~\ref{prop:delay} adds is that the remainder is bounded \emph{uniformly in $\delta$}, although the number of steps and the constants in the classical expansion grow without bound as $\delta\to0$. This uniformity is what allows conclusions about the limit $\delta\to0$.

For explicit Euler, $C=-\frac12f_qf_q'$ and $g=f_q'/(2f_q)$, so that $D_n\approx\frac h2\ln\bigl(f_q(y_n)/f_q(\delta)\bigr)$ and $\kappa=q/2$. For a general first-order method $\kappa$ may vanish; the delay is then $O(h)$ uniformly in $\delta$, so that $h\ln(1/\delta)$ is an upper bound for its growth, attained by the Euler methods.

The interpretation is simple. The delay is accumulated where $g$ is large. For first-order methods $g\sim\kappa/y$ and the contribution of the long initial phase is only logarithmic, because there the solution varies on the slow time scale $1/y$ and the local errors are correspondingly tiny. For $p\ge2$ the density is integrable at $y=0$ and the delay is produced in the transition itself, in a $\delta$-independent window. In neither case is there an accumulation proportional to the length $1/\delta$ of the unstable interval.

\begin{remark}[The choice of variable]\label{rem:variables}
The logarithmic growth for $p=1$ is a property of the variable in which the method is applied, not of the problem. Let explicit Euler be applied to $w'=F(w)$, where $w=\varphi(y)$ and $F(\varphi(y))=\varphi'(y)f_q(y)$. The one-step time defect of Step~1 is then $\tau=-\frac{h^2}{2}F'(w)+O(h^3)$, and the computation leading to \eqref{eq:delay} gives
\[
D_n\approx\frac h2\ln\frac{\varphi'(y_n)\,f_q(y_n)}{\varphi'(\delta)\,f_q(\delta)} .
\]
For $w=y$ this is the formula above, which grows like $\frac q2h\ln(1/\delta)$. For $w=\ln y$ the growth is $\frac{q-1}2h\ln(1/\delta)$. For $w=y^{1-q}$, which turns \eqref{eq:fq} into $w'=-(q-1)(1-y)$, the delay $D_n\approx\frac h2\ln\bigl((1-y_n)/(1-\delta)\bigr)$ does not grow at all: in this variable the slow phase is the nearly linear motion $w\approx\delta^{1-q}-(q-1)t$, the regularized blow-up of Section~\ref{sec:structure}, which Euler integrates almost exactly. For $q=2$, $h=0.25$ and $\delta=10^{-2},10^{-4},10^{-6}$ the computed delays at the first step with $y_k\ge1/2$ are $0.887$, $2.036$, $3.194$ in the variable $y$, $0.403$, $0.977$, $1.553$ in $\ln y$, and $-0.093$, $-0.090$, $-0.086$ in $1/y$. Over $h=0.25,0.125$, $q=2$ ($10^{-6}\le\delta\le10^{-2}$) and $q=3$ ($10^{-4}\le\delta\le10^{-2}$) the computed delays are within $1.2\%$ of the formula above in all three variables, and the scaled remainder $(D_n-\text{leading term})/h^2$ is at most $0.11$ in absolute value and does not grow as $\delta$ decreases. We have not proved this uniformity for the transformed equations. For $p\ge2$ the choice of variable matters less, since the delay is then bounded uniformly in $\delta$ already in the variable $y$.
\end{remark}

Table~\ref{tab:g} lists the densities, computed symbolically.

\begin{table}[!tbp]
\centering\small
\caption{Delay densities $g=-C/f^2$ and antiderivatives $G$ ($D\approx h^p[G(y)-G(\delta)]$) for $f=y^2-y^3$.}
\label{tab:g}
\renewcommand{\arraystretch}{1.35}
\begin{tabular}{lcccl}
\toprule
method & $p$ & $g(y)$ & $\kappa$ & $G(y)$\\
\midrule
explicit Euler & 1 & $\dfrac{2-3y}{2y(1-y)}$ & $1$ & $\ln y+\tfrac12\ln(1-y)$\\
implicit Euler & 1 & $-\dfrac{2-3y}{2y(1-y)}$ & $-1$ & $-\ln y-\tfrac12\ln(1-y)$\\
Heun & 2 & $\dfrac{6y^2-8y+3}{6(1-y)}$ & -- & $-\tfrac{y^2}{2}+\tfrac{y}{3}-\tfrac16\ln(1-y)$\\
RK4 & 4 & $\dfrac{y^2(15y^4-40y^3+55y^2-38y+10)}{240(1-y)}$ & -- & see code\\
\bottomrule
\end{tabular}
\end{table}

\begin{table}[!tbp]
\centering\small
\caption{Exact phase delay $D_k$ at the first step with $y_k\ge1/2$, $h=0.25$, and ratio to the leading term $h^p\int_\delta^{y_k}g$ of \eqref{eq:delay}. Iterates accumulated with compensated summation, $D_k$ accumulated step by step (see text).}
\label{tab:phase}
\begin{tabular}{lcccccc}
\toprule
 & \multicolumn{2}{c}{$\delta=10^{-2}$} & \multicolumn{2}{c}{$\delta=10^{-4}$} & \multicolumn{2}{c}{$\delta=10^{-6}$}\\
\cmidrule(lr){2-3}\cmidrule(lr){4-5}\cmidrule(lr){6-7}
method & $D$ & ratio & $D$ & ratio & $D$ & ratio\\
\midrule
explicit Euler & $8.868\cdot10^{-1}$ & 0.993 & $2.036\cdot10^{0}$ & 0.997 & $3.194\cdot10^{0}$ & 0.998 \\
implicit Euler & $-9.011\cdot10^{-1}$ & 1.008 & $-2.053\cdot10^{0}$ & 1.003 & $-3.205\cdot10^{0}$ & 1.002 \\
Heun & $9.316\cdot10^{-3}$ & 0.960 & $9.483\cdot10^{-3}$ & 0.962 & $9.661\cdot10^{-3}$ & 0.961 \\
RK4 & $1.964\cdot10^{-6}$ & 0.968 & $1.883\cdot10^{-6}$ & 0.966 & $1.985\cdot10^{-6}$ & 0.968 \\
 
\bottomrule
\end{tabular}
\end{table}

\subsection{Numerical confirmation}

\begin{figure}[!tbp]
\centering
\includegraphics[width=\textwidth]{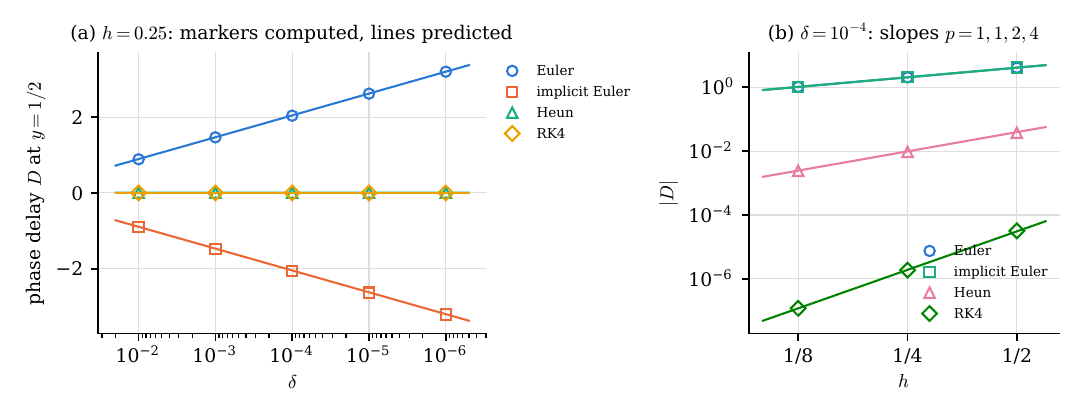}
\caption{(a) Phase delay at $y=1/2$ versus $\delta$ for $h=0.25$ (markers: computed exactly via Lemma~\ref{lem:phase}; lines: \eqref{eq:delay}). (b) Delay versus $h$ for $\delta=10^{-4}$.}
\label{fig:phase}
\end{figure}

Table~\ref{tab:phase} and Fig.~\ref{fig:phase} compare the exact delays with the leading term of \eqref{eq:delay}. The agreement is within $0.8\%$ for the Euler methods and within $4\%$ for Heun and RK4 at $h=0.25$, improving as $h\to0$ with the predicted slopes $p=1,2,4$ (Fig.~\ref{fig:phase}b). For first-order methods the delay grows by $\kappa h\ln10$ per decade of $\delta$; for Heun and RK4 it does not grow (it is constant to leading order).

Table~\ref{tab:remainder} tests the statement of Theorem~\ref{prop:delay} itself: the scaled remainder $\bigl(D_k-h^p\int_\delta^{y_k}g\bigr)/h^{p+1}$ remains of comparable magnitude over four decades of $\delta$ ($q=2$), two decades ($q=3$) and both step sizes, for every method, consistent with the uniform bound; the three entries marked $^\dagger$ are affected by rounding, as discussed below.

These computations require some care with rounding. By Proposition~\ref{prop:sens} (and its analogue for $q\ge2$), a relative rounding error $u$ at level $y$ shifts the phase by about $u/y^{q-1}$. Over the $10^7$--$10^8$ steps of the slow phase these shifts add up to an accumulated rounding delay comparable with the truncation delay of RK4. We therefore accumulate the iterates with compensated (Kahan) summation of the increments $\Phi_h(y_k)-y_k$, which are computed directly, and accumulate $D_k$ step by step as $D_{k+1}=D_k+h-[\Theta_q(y_{k+1})-\Theta_q(y_k)]$, with the difference of $\Theta_q$ evaluated without cancellation through \texttt{log1p} and \texttt{expm1}. With plain summation, for example, RK4 with $h=0.125$ and $\delta=10^{-6}$ gives half the correct delay; with compensated summation the ratio to the prediction is $0.983$. The remaining visible rounding effects are the RK4 entries for $q=3$ marked $^\dagger$ in Table~\ref{tab:remainder}, where the delay itself is below $3\cdot10^{-7}$ ($h=0.25$) and $2\cdot10^{-8}$ ($h=0.125$), and the remainder term is a few times $10^{-9}$.

\begin{table}[!tbp]
\centering\small
\caption{Scaled remainder $\bigl(D_k-h^p\int_\delta^{y_k}g\bigr)/h^{p+1}$ of Theorem~\ref{prop:delay} at the first step with $y_k\ge1/2$, for $f_q=y^q(1-y)$. The theorem bounds it uniformly in $\delta$; the table tests that it does not grow as $\delta$ decreases. $^\dagger$: affected by rounding.}
\label{tab:remainder}
\begin{tabular}{llccccc}
\toprule
$q=2$ & $h$ & $10^{-2}$ & $10^{-3}$ & $10^{-4}$ & $10^{-5}$ & $10^{-6}$ \\
\midrule
explicit Euler & 0.25 & $-0.102$ & $-0.106$ & $-0.107$ & $-0.106$ & $-0.105$ \\
 & 0.125 & $-0.103$ & $-0.107$ & $-0.108$ & $-0.108$ & $-0.108$ \\
implicit Euler & 0.25 & $-0.108$ & $-0.112$ & $-0.112$ & $-0.112$ & $-0.112$ \\
 & 0.125 & $-0.106$ & $-0.110$ & $-0.111$ & $-0.111$ & $-0.111$ \\
Heun & 0.25 & $-0.025$ & $-0.025$ & $-0.024$ & $-0.025$ & $-0.025$ \\
 & 0.125 & $-0.024$ & $-0.025$ & $-0.025$ & $-0.025$ & $-0.025$ \\
RK4 & 0.25 & $-6.73\cdot10^{-5}$ & $-6.79\cdot10^{-5}$ & $-6.83\cdot10^{-5}$ & $-6.60\cdot10^{-5}$ & $-6.68\cdot10^{-5}$ \\
 & 0.125 & $-6.85\cdot10^{-5}$ & $-6.80\cdot10^{-5}$ & $-6.84\cdot10^{-5}$ & $-6.77\cdot10^{-5}$ & $-6.65\cdot10^{-5}$ \\
\midrule
$q=3$ & $h$ & $10^{-2}$ & $3\cdot10^{-3}$ & $10^{-3}$ & $3\cdot10^{-4}$ & $10^{-4}$ \\
\midrule
explicit Euler & 0.25 & $-0.063$ & $-0.064$ & $-0.064$ & $-0.064$ & $-0.064$ \\
 & 0.125 & $-0.064$ & $-0.064$ & $-0.064$ & $-0.064$ & $-0.065$ \\
implicit Euler & 0.25 & $-0.066$ & $-0.067$ & $-0.067$ & $-0.067$ & $-0.067$ \\
 & 0.125 & $-0.066$ & $-0.066$ & $-0.066$ & $-0.066$ & $-0.066$ \\
Heun & 0.25 & $-0.006$ & $-0.006$ & $-0.006$ & $-0.006$ & $-0.006$ \\
 & 0.125 & $-0.006$ & $-0.006$ & $-0.006$ & $-0.006$ & $-0.006$ \\
RK4 & 0.25 & $-5.73\cdot10^{-6}$ & $-5.27\cdot10^{-6}$ & $-5.42\cdot10^{-6}$ & $-5.86\cdot10^{-6}$ & $-2.46\cdot10^{-6}$$^\dagger$ \\
 & 0.125 & $-5.01\cdot10^{-6}$ & $-5.02\cdot10^{-6}$ & $-4.35\cdot10^{-6}$ & $4.77\cdot10^{-6}$$^\dagger$ & $8.24\cdot10^{-5}$$^\dagger$ \\
 
\bottomrule
\end{tabular}
\end{table}

%=====================================================================
\section{The MATLAB ODE suite on the flame problem: cost and accuracy}\label{sec:adaptive}

We repeated the experiment of \cite{Moler2014} with MATLAB R2023b (Update 3): problem \eqref{eq:flame} with $\mathrm{RelTol}=10^{-4}$ and the default $\mathrm{AbsTol}=10^{-6}$, for the eight solvers \texttt{ode45}, \texttt{ode23}, \texttt{ode113}, \texttt{ode23s}, \texttt{ode15s} (default NDF formulas and with the option \texttt{'BDF','on'}), \texttt{ode23t} and \texttt{ode23tb} \cite{ShampineReichelt1997}. We call this the benchmark setting, and $\mathrm{RelTol}=10^{-4}$, $\mathrm{AbsTol}=10^{-6}$ the benchmark tolerance. The step statistics reproduce the numbers reported in \cite{Moler2014} exactly (Table~\ref{tab:moler1e4}): \texttt{ode45} takes $3040$ steps and $20179$ function evaluations, \texttt{ode23s} $99$ steps and $412$ evaluations, and \texttt{ode15s} $140$ steps, $39$ failed attempts and $347$ evaluations. For every accepted step we also record whether it lies in the slow phase ($t<c(\delta)-20$) or in the stiff phase ($t>c(\delta)+20$), the error against the exact solution, and the exact delay $D_k$ of Lemma~\ref{lem:phase}. All solvers were used with their default options except the tolerances. Each MATLAB solver returns a solution structure; step statistics and nodal errors use its accepted steps (\texttt{sol.x}; the \texttt{Refine} option applies only to array output), and the additional output times of Section~\ref{sec:sundman} are evaluated with \texttt{deval}. Work is measured by the solvers' own counts of accepted steps, failed steps and function evaluations. SciPy's \texttt{solve\_ivp} likewise returns the accepted steps when no output times are prescribed. The scripts are in the code repository \cite{BorosCode2026}. The SciPy (version 1.18) counterparts RK45, RK23, Radau and BDF give the same picture.

\begin{table}[t]
\centering\small
\setlength{\tabcolsep}{4pt}
\caption{The benchmark run ($\delta=10^{-4}$, $\mathrm{RelTol}=10^{-4}$, $\mathrm{AbsTol}=10^{-6}$), MATLAB R2023b: accepted and failed steps, function evaluations, Jacobian evaluations, LU decompositions, steps in the slow and in the stiff phase, front delay $|D|$ at $y=1/2$ and maximal error. ``(warm)'': started from the exact solution at $y=0.01$.}
\label{tab:moler1e4}
\begin{tabular}{lrrrrrrrcc}
\toprule
solver & steps & failed & $f$-evals & Jac. & LU & slow & stiff & $|D|$ & error\\
\midrule
\texttt{ode45} & 3040 & 323 & 20179 & 0 & 0 & 17 & 3007 & $2.2\cdot10^{-1}$ & $3.0\cdot10^{-2}$ \\
\texttt{ode23} & 4052 & 12 & 12193 & 0 & 0 & 39 & 3992 & $2.8\cdot10^{1}$ & $9.5\cdot10^{-1}$ \\
\texttt{ode113} & 6239 & 823 & 13302 & 0 & 0 & 28 & 6202 & $3.3\cdot10^{1}$ & $9.5\cdot10^{-1}$ \\
\texttt{ode23s} & 99 & 7 & 412 & 99 & 106 & 87 & 8 & $2.9\cdot10^{1}$ & $9.5\cdot10^{-1}$ \\
\texttt{ode15s} & 140 & 39 & 347 & 2 & 53 & 132 & 8 & $3.3\cdot10^{2}$ & $1.0\cdot10^{0}$ \\
\texttt{ode15s} (BDF) & 167 & 40 & 400 & 2 & 53 & 159 & 8 & $4.9\cdot10^{2}$ & $1.0\cdot10^{0}$ \\
\texttt{ode23t} & 157 & 4 & 337 & 2 & 54 & 150 & 7 & $2.2\cdot10^{2}$ & $9.9\cdot10^{-1}$ \\
\texttt{ode23tb} & 134 & 9 & 551 & 2 & 47 & 127 & 7 & $7.5\cdot10^{1}$ & $9.8\cdot10^{-1}$ \\
\midrule
\texttt{ode45} (warm) & 3028 & 319 & 20083 & 0 & 0 & 5 & 3007 & $4.1\cdot10^{-4}$ & $1.1\cdot10^{-4}$ \\
\texttt{ode23s} (warm) & 80 & 0 & 322 & 80 & 80 & 17 & 12 & $6.3\cdot10^{-2}$ & $9.5\cdot10^{-3}$ \\
\texttt{ode15s} (warm) & 97 & 15 & 205 & 2 & 28 & 24 & 21 & $1.4\cdot10^{-1}$ & $2.0\cdot10^{-2}$ \\
 
\bottomrule
\end{tabular}
\end{table}

\subsection{Cost}

\paragraph{In what sense the problem is stiff} The problem is stiff only after the front, and only for small $\delta$. Before the front $f'(y)=2y-3y^2$ is positive, and of order $\delta$ in the initial part of the slow phase: the problem is mildly unstable, not stiff. After the front the solution is constant to within the tolerance, while nearby solutions approach it like $e^{f'(1)t}=e^{-t}$. The ratio of the length $L_{\mathrm{stiff}}\approx1/\delta$ of the remaining interval to this time scale, $|f'(1)|\,L_{\mathrm{stiff}}\approx1/\delta$, measures the stiffness. The operational criterion confirms it: an explicit code whose step is limited by stability does not react to the tolerance. For $\delta=10^{-4}$, \texttt{ode45} takes between $3005$ and $3013$ steps in the stiff phase for every tolerance from $10^{-3}$ to $10^{-12}$ (standard strategy of Section~\ref{sec:wp}), and \texttt{ode113} between $5952$ and $6271$, while the implicit codes take $5$--$54$. As $\delta$ decreases from $2\cdot10^{-2}$ to $10^{-6}$, the stiff-phase steps of \texttt{ode45} at the benchmark tolerance grow from $9$ to $301416$, by a factor of about $10$ per decade of $\delta$ for $\delta\le10^{-3}$, while those of \texttt{ode15s} stay between $5$ and $11$. For $\delta=10^{-2}$ the problem is thus hardly stiff ($24$ steps of \texttt{ode45} against $9$ of \texttt{ode15s} in this phase).

\paragraph{The stiff phase} After the front, $y\approx1$ and $f'(1)=-1$. The solution is constant to within the tolerance, so an explicit method is limited by absolute stability, $h\le\beta$ (a mechanism well known from the analysis of step size control on stiff problems \cite{Shampine1977,HallHigham1988}, \cite[Sect.~IV.2]{HairerWanner1996}), where $\beta$ is the real stability boundary of the method (the largest $x$ with $|R(-x)|\le1$ for its stability function $R$). From the Butcher tableaux, $\beta=3.3066$ for the Dormand--Prince 5(4) pair of \texttt{ode45} (which propagates the fifth-order solution) and $\beta=2.5127$ for the Bogacki--Shampine 3(2) pair of \texttt{ode23}. The stiff phase has length
\[
L_{\mathrm{stiff}}=\frac2\delta-c(\delta)=\frac1\delta-\ln\frac1\delta+2+O(\delta),
\]
so an explicit code needs about $L_{\mathrm{stiff}}/\beta$ steps there. Table~\ref{tab:stiffsteps} shows that this prediction is correct to within $2\%$ for \texttt{ode45} ($0.5\%$ for $\delta\le10^{-4}$) and $1.8\%$ for \texttt{ode23}. For $\delta=10^{-4}$ it gives $3022$; \texttt{ode45} takes $3007$ of its $3040$ steps in the stiff phase. The step count of \texttt{ode45} is thus the length of the stiff phase divided by the stability boundary of the Dormand--Prince pair; for this problem the length $L_{\mathrm{stiff}}$ is known in closed form, so the count can be predicted rather than only explained. The step-size controller does not settle exactly at $\beta$: a step slightly larger than $\beta$ amplifies the error, the error estimate grows, and the next step is reduced or rejected. Whether the controller has a stable equilibrium near $\beta$ depends on the method and the controller (step-control stability \cite{HallHigham1988}). The rejected steps of \texttt{ode113} hardly depend on the tolerance ($819$--$952$ for $\mathrm{tol}=10^{-3},\dots,10^{-12}$, $\delta=10^{-4}$), nor do those of \texttt{ode45} for $\mathrm{tol}\le10^{-5}$ ($161$--$214$), which indicates that they arise from this interaction with the stability boundary. They affect only the cost; the phase error is created before the front, in the slow phase. The variable-order Adams code \texttt{ode113} has no single stability boundary; it advances on average only $1.6$--$1.7$ time units per step in the stiff phase.

\paragraph{The slow phase} Near $y=0$ the equation is invariant under $y\mapsto\lambda y$, $t\mapsto t/\lambda$, and so is relative error control. The accepted steps therefore satisfy $h_ky_k\approx\mathrm{const}$, the ratio $y_{k+1}/y_k$ is roughly constant, and the slow phase of length $\approx1/\delta$ costs only $O(\ln(1/\delta))$ steps. In the benchmark setting $\mathrm{AbsTol}/\mathrm{RelTol}=10^{-2}$, so absolute error control takes over for $y<10^{-2}$. The explicit codes need $3$--$47$ slow-phase steps for $\delta$ from $2\cdot10^{-2}$ to $10^{-6}$, with step sizes beyond $10^3$ (Fig.~\ref{fig:moler}a).

\paragraph{Stiff solvers} Implicit solvers are not restricted in the stiff phase: their steps grow geometrically, and the phase costs them $3$--$13$ steps. Their cost is the slow phase. Of the $99$ steps of \texttt{ode23s} at $\delta=10^{-4}$, $87$ lie before the front, and of the $140$ steps of \texttt{ode15s}, $132$. The totals grow only logarithmically, from $64$ to $118$ steps for \texttt{ode23s} and from $75$ to $147$ for \texttt{ode15s} as $\delta$ decreases from $2\cdot10^{-2}$ to $10^{-6}$ (Fig.~\ref{fig:moler}b).

\begin{table}[t]
\centering\small
\caption{Stiff-phase steps of the explicit MATLAB solvers at $\mathrm{RelTol}=10^{-4}$, $\mathrm{AbsTol}=10^{-6}$, and the prediction $L_{\mathrm{stiff}}/\beta$.}
\label{tab:stiffsteps}
\begin{tabular}{cccccccc}
\toprule
 & & \multicolumn{2}{c}{\texttt{ode23} ($\beta=2.513$)} & \multicolumn{2}{c}{\texttt{ode45} ($\beta=3.307$)} & \texttt{ode113}\\
\cmidrule(lr){3-4}\cmidrule(lr){5-6}\cmidrule(lr){7-7}
$\delta$ & $L_{\mathrm{stiff}}$ & steps & predicted & steps & predicted & steps\\
\midrule
$10^{-3}$ & 995 & 389 & 396 & 295 & 301 & 583 \\
$10^{-4}$ & 9993 & 3992 & 3977 & 3007 & 3022 & 6202 \\
$10^{-5}$ & 99990 & 39581 & 39794 & 30135 & 30240 & 60256 \\
$10^{-6}$ & 999988 & 391721 & 397974 & 301416 & 302422 & 583383 \\
 
\bottomrule
\end{tabular}
\end{table}

\subsection{Accuracy}

Standard adaptive codes control the error per step, and their step size follows the local time scale. The following identity, which extends Step~2 of the proof of Theorem~\ref{prop:delay} to arbitrary step sizes, shows what this does to the phase.

\begin{proposition}[Delay of an arbitrary step sequence]\label{prop:defects}
Let $(t_k,y_k)$, $k=0,\dots,n$, be any sequence of points with $t_0=0$, $y_0=\delta$, step sizes $h_k=t_{k+1}-t_k>0$ and $\delta\le y_k<1$ for all $k$, computed by any method. Let $\ell_k:=y_{k+1}-\mathcal{F}_{h_k}(y_k)$ be the defect of step $k$ with respect to the exact flow $\mathcal{F}_h$ of \eqref{eq:fq} from the current computed state. Then
\[
D_n=-\sum_{k<n}\tau_k,\qquad \tau_k:=\int_{\mathcal{F}_{h_k}(y_k)}^{y_{k+1}}\frac{\dd x}{f_q(x)}=\frac{\ell_k}{f_q(\xi_k)},
\]
with $\xi_k$ between $\mathcal{F}_{h_k}(y_k)$ and $y_{k+1}$.
\end{proposition}

\begin{proof}
Since $\mathcal{F}_{h}(y)\in[y,1)$ and $T_\delta(\mathcal{F}_h(y))=T_\delta(y)+h$, we have $D_{k+1}-D_k=h_k-\bigl[T_\delta(y_{k+1})-T_\delta(y_k)\bigr]=-\bigl[T_\delta(y_{k+1})-T_\delta(\mathcal{F}_{h_k}(y_k))\bigr]=-\tau_k$, and $D_0=0$. The last expression for $\tau_k$ is the mean value theorem.
\end{proof}

For a one-step method $\ell_k$ is the local error. For a multistep method it differs from the local truncation error, and for any adaptive code it differs from the internally estimated error, so an error test controls $\ell_k$ only approximately. The delay is thus exactly the sum of the defects, each converted into time by the factor $1/f_q$ at the place where it is committed. An existing time shift is preserved by the exact subsequent evolution, since every later solution is a translate of the same profile; later numerical defects can add to it or cancel it, but the flow supplies no correction. If a code meets its error test in the true local error, $|\ell_k|\le\max(\mathrm{RelTol}\cdot y_k,\mathrm{AbsTol})$, the admissible time defect in the slow phase is about $\mathrm{RelTol}\,y_k^{1-q}$ or $\mathrm{AbsTol}\,y_k^{-q}$ per step: it is largest at the start, where $y_k\approx\delta$. The simplest adaptive method shows in closed form what this does.

\begin{figure}[t]
\centering
\includegraphics[width=\textwidth]{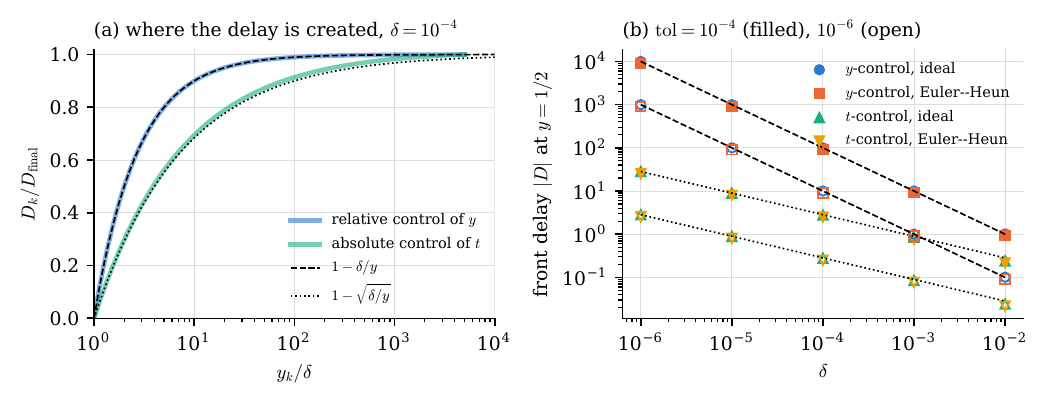}
\caption{Adaptive explicit Euler for \eqref{eq:flame}, with relative control of $y$ (local error $\mathrm{tol}\cdot y$) and with time as a dependent variable and absolute control of $t$ (Section~\ref{sec:sundman}; local error of $t$ equal to $\mathrm{tol}$). ``Ideal'': the leading local error term equals the tolerance; ``Euler--Heun'': Heun's method as error estimator with a standard step-size controller. (a) Fraction of the final delay accumulated by the level $y_k$, $\mathrm{tol}=10^{-4}$. (b) Front delay; lines: $\sqrt{\mathrm{tol}}/\delta$ (dashed) and $\sqrt{8\,\mathrm{tol}/\delta}$ (dotted).}
\label{fig:euler}
\end{figure}

\paragraph{An adaptive Euler method} Let explicit Euler choose each step so that its leading local error $\frac12h_k^2f'(y_k)f(y_k)$ equals $\mathrm{tol}\cdot y_k$. In the slow phase $f(y)\approx y^2$ and $f'(y)\approx2y$, so $h_k\approx\sqrt{\mathrm{tol}}/y_k$ and $y_{k+1}=y_k+h_kf(y_k)\approx(1+\sqrt{\mathrm{tol}})\,y_k$: the iterates grow geometrically, and the slow phase costs $O(\ln(1/\delta)/\sqrt{\mathrm{tol}})$ steps. By Proposition~\ref{prop:sens}, each local error $\mathrm{tol}\cdot y_k$ is a time defect $\mathrm{tol}\cdot y_k/f(y_k)\approx\mathrm{tol}/y_k$, and the defects add up (Proposition~\ref{prop:defects}):
\[
D\approx\sum_{k\ge0}\frac{\mathrm{tol}}{\delta\,(1+\sqrt{\mathrm{tol}})^{k}}\approx\frac{\sqrt{\mathrm{tol}}}{\delta}.
\]
The sum is dominated by its first terms: the part of $D$ accumulated before the solution reaches the level $y$ is the fraction $1-\delta/y$ of the total, so half of the delay is created before $y=2\delta$ and $90\%$ before $y=10\delta$ (Fig.~\ref{fig:euler}a). For $\mathrm{tol}=10^{-4}$ and $\delta=10^{-4}$ the front comes $100$ time units late. The exact delays agree with $\sqrt{\mathrm{tol}}/\delta$ to within $1.2\%$ for $\mathrm{tol}=10^{-4},10^{-6}$ and $\delta=10^{-2},\dots,10^{-6}$; with a practical controller (Euler propagated, Heun as error estimator, rejected steps repeated) they are $0.90$--$0.92$ times this value (Fig.~\ref{fig:euler}b).

\begin{figure}[t]
\centering
\includegraphics[width=\textwidth]{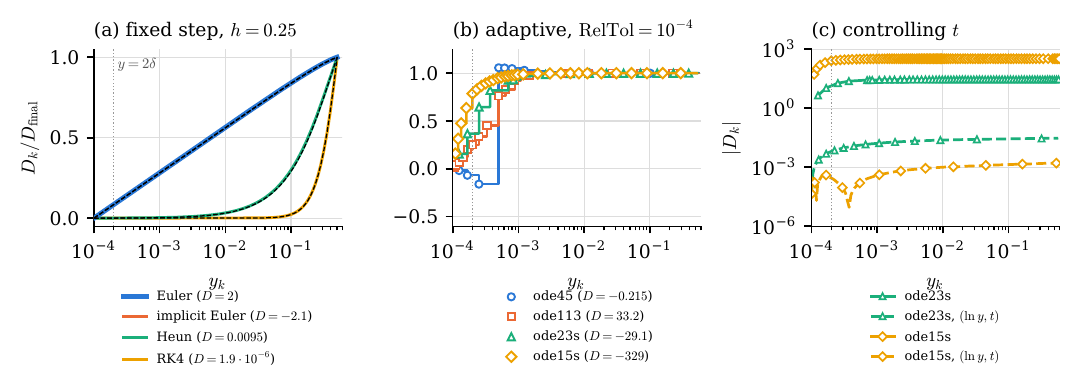}
\caption{How the phase delay accumulates, $\delta=10^{-4}$; the dotted line marks $y=2\delta$. (a) Fixed step $h=0.25$: $D_k/D_{\mathrm{final}}$ (solid) and the leading term of \eqref{eq:delay} (dashed); normalized, implicit Euler coincides with explicit Euler. (b) MATLAB R2023b, $\mathrm{RelTol}=10^{-4}$, $\mathrm{AbsTol}=10^{-6}$: $D_k/D_{\mathrm{final}}$ at the accepted steps. (c) $|D_k|$ for \texttt{ode23s} and \texttt{ode15s}, standard formulation (solid) and time as a dependent variable \eqref{eq:sundmanlog} with $\mathrm{AbsTol}=10^{-4}$ for $\ln y$ and $t$ (dashed). $D_{\mathrm{final}}$: delay at the first step with $y_k\ge1/2$.}
\label{fig:profile}
\end{figure}

\paragraph{A heuristic delay law} For a code that controls an error estimate of order $p+1$, the same computation gives a constant growth factor $1+(\mathrm{tol}/C)^{1/(p+1)}$ per step and hence a delay proportional to $\mathrm{tol}^{p/(p+1)}/\delta$. In the absolute regime the admissible local error is $\mathrm{AbsTol}$, the time defect $\mathrm{AbsTol}/y^2$, and the growth factor per step depends on $y$; summing the defects in the same way gives
\begin{equation}\label{eq:adaptive}
D\approx K\,\frac{\mathrm{RelTol}^{\alpha}}{\delta}\quad\text{(relative regime)},\qquad D\approx K'\,\frac{\mathrm{AbsTol}^{\alpha}}{\delta^{1+\alpha}}\quad\text{(when }\mathrm{AbsTol}>\mathrm{RelTol}\cdot\delta\text{)},
\end{equation}
where $\alpha=p/(p+1)$ and $K,K'$ are method-dependent constants. If the code propagates the more accurate solution of an embedded pair (local extrapolation), the actual local error is one order smaller than the controlled one, and the same computation gives $\alpha=1$. Equation \eqref{eq:adaptive} is a model, not a theorem: in the relative regime it assumes that every step in the slow phase commits a local error equal to the tolerance and that the ratio $y_{k+1}/y_k$ is constant, so that Proposition~\ref{prop:defects} turns the delay into a geometric sum; in the absolute regime the growth factor per step depends on $y$, and the sum is evaluated with that state-dependent factor. The exponents $\alpha$ are those of the classical analysis of tolerance proportionality \cite{Higham1993,SoderlindWang2006}; the factor $1/\delta$ is specific to the sensitivity of the front. The MATLAB runs agree with these exponents (all exponents quoted are least-squares slopes in log--log scale; the fit ranges are listed in the README file of the code repository \cite{BorosCode2026}). The delays of the second-order codes \texttt{ode23s}, \texttt{ode23t} and \texttt{ode23tb} scale like $\mathrm{RelTol}^{0.67}$ and $\mathrm{AbsTol}^{0.67}$, those of \texttt{ode15s} (orders up to $5$) like $\mathrm{RelTol}^{0.82\text{--}0.84}$ and $\mathrm{AbsTol}^{0.86}$, and those of \texttt{ode23} and \texttt{ode113}, which propagate the more accurate solution, like $\mathrm{RelTol}^{0.97\text{--}1.02}$. For \texttt{ode45}, which also propagates the more accurate solution, the exponents scatter between $0.81$ and $1.08$. In the absolute regime the delays of the three second-order codes grow by factors of $2130$--$2290$ (medians) per two decades of $\delta$, against $100^{5/3}\approx2150$ predicted. As in the Euler example, the delay is created at the very beginning: for SciPy's RK45 with $\mathrm{RelTol}=10^{-12}$ and $\delta=10^{-4}$, $50\%$ of the final delay has been incurred by the time $y=2\delta$, and $90\%$ by $y=10\delta$. The flow does not damp this delay, since every solution is a translate of the same profile, and the later steps change it little (Fig.~\ref{fig:profile}b).

\begin{figure}[t]
\centering
\includegraphics[width=\textwidth]{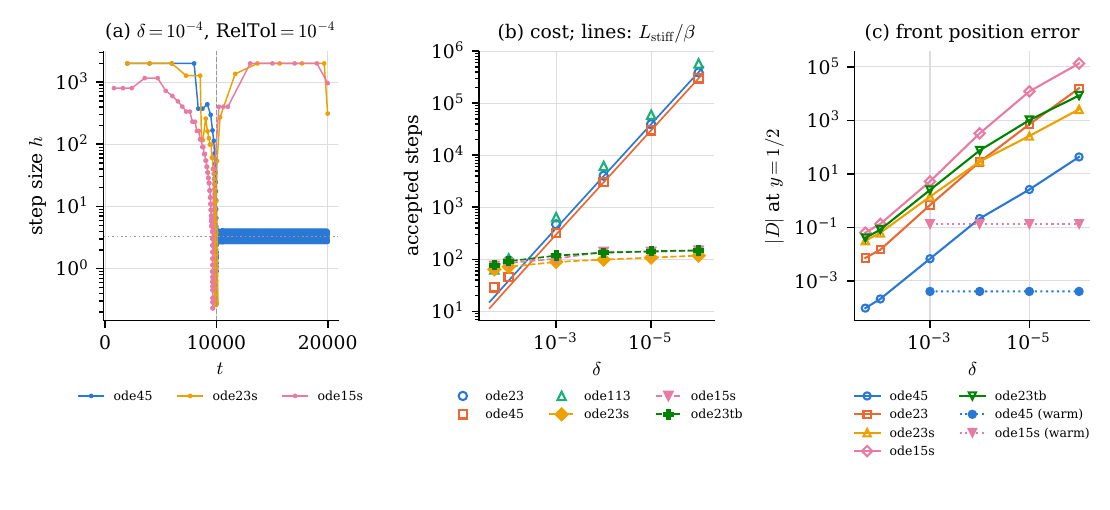}
\caption{Benchmark setting, MATLAB R2023b, $\mathrm{RelTol}=10^{-4}$, $\mathrm{AbsTol}=10^{-6}$. (a) Step sizes for $\delta=10^{-4}$; the dashed line marks the front, the dotted line the stability boundary of the Dormand--Prince pair. (b) Accepted steps versus $\delta$; lines: $L_{\mathrm{stiff}}/\beta$ for \texttt{ode45} and \texttt{ode23}. (c) Front delay at $y=1/2$; ``warm'': started from the exact solution at $y=0.01$.}
\label{fig:moler}
\end{figure}

\paragraph{The benchmark tolerance} Tables~\ref{tab:moler1e4} and~\ref{tab:sweep} and Fig.~\ref{fig:moler}(c) give the delays at the benchmark tolerance. For $\delta=10^{-4}$ only \texttt{ode45} is accurate: it places the front $0.22$ time units early, with maximal error $0.030$. The solutions that the stiff solvers compute so cheaply are not. \texttt{ode23s} places the front $29$ time units early and \texttt{ode15s} $329$ time units early, with maximal errors $0.95$ and $1.00$. The other solvers misplace it by $28$ (\texttt{ode23}) to $487$ (\texttt{ode15s} with BDF) time units. On a plot of the whole interval $[0,2\cdot10^4]$ a shift of $329$ is $1.6\%$ of the width and hardly visible, so the stiff solutions look correct. The delays of \texttt{ode45}, \texttt{ode23s} and \texttt{ode15s} grow by factors between $9$ and $62$ per decade of $\delta$, between the two regimes of \eqref{eq:adaptive}. For $\delta=10^{-6}$ every solver, including \texttt{ode45}, has an error close to $1$. A comparison of step counts at fixed tolerance thus compares solutions of very different quality.

The delays do not depend on the level at which the front is measured: for the eight solvers at $\delta=10^{-4}$ the delays at the first computed points with $y_k\ge0.1$, $0.5$ and $0.9$ differ by at most $0.1\%$. Nor does it help to locate the front with an event function, since the event function sees the delayed computed solution. With the event location of SciPy for $y=1/2$ (RK45, Radau, BDF and LSODA at the benchmark tolerance) the located times differ from the exact one by $-0.11$, $0.09$, $-131.5$ and $-148.2$, equal to $D_{\mathrm{front}}$ of the same runs to within $0.2\%$.

\paragraph{The sign of the delay} The direction of the error is also set in the slow phase. There $f'(y)>0$ and the solution is convex; in fact all its derivatives are positive, since $y'\approx y^2$ gives $y^{(k)}\approx k!\,y^{k+1}$. By Assumption~(S) the leading local error of a one-step method near $y=0$ is $E(0,0)\,y\,z^{p+1}$ with a constant $E(0,0)$ that depends only on the method, so all local errors in the slow phase have one sign and the delays add up instead of cancelling. Explicit Euler falls short of the convex solution (local error $-\frac12h^2y''$) and its front is late; implicit Euler overshoots ($+\frac12h^2y''$) and its front is early ($\kappa=\pm1$ in Table~\ref{tab:g}). The MATLAB codes follow the same pattern. Over all completed runs without warm start, the implicit codes \texttt{ode15s} (NDF and BDF), \texttt{ode23t} and \texttt{ode23tb} place the front early in every run ($33$, $33$, $31$ and $33$ runs), and \texttt{ode23s} in $31$ of $33$ runs; in the two exceptions ($q=3$, $\mathrm{RelTol}=10^{-4}$, $\delta\le10^{-4}$) it misses the front altogether, with errors above $0.99$. The explicit codes \texttt{ode23} and \texttt{ode113} place the front late in all $27$ runs. \texttt{ode45}, which propagates the more accurate solution of its pair, has both signs ($15$ early, $12$ late).

\paragraph{Where the delay is created} Figure~\ref{fig:profile} follows $D_k$ along the trajectory for $\delta=10^{-4}$. With fixed steps (a) the delay follows the leading term of Theorem~\ref{prop:delay}: for the Euler methods it grows linearly in $\ln y_k$ over the whole slow phase, for Heun and RK4 it is created in the transition, $y\gtrsim0.1$. The adaptive codes at the benchmark tolerance (b) behave like the adaptive Euler method. By the time $y$ has reached $10\delta$, after $6$ to $26$ accepted steps, every code except \texttt{ode45} has incurred $92$--$99\%$ of its final delay, and afterwards $D_k$ changes by at most $5\%$ of the final value. For \texttt{ode45} the final delay $-0.22$ is the small remainder of contributions of both signs. With time as a dependent variable and absolute control of $t$ (c), the defects are spread over the slow phase instead, and $|D_k|$ stays below $0.03$ (\texttt{ode23s}) and $0.002$ (\texttt{ode15s}) up to $y=1/2$, against $29$ and $329$ in the standard formulation.

\begin{table}[t]
\centering\small
\setlength{\tabcolsep}{4pt}
\caption{MATLAB solvers at the benchmark tolerance for varying $\delta$: accepted steps, front delay $|D|$ at $y=1/2$ and maximal error.}
\label{tab:sweep}
\begin{tabular}{cccccccccc}
\toprule
 & \multicolumn{3}{c}{\texttt{ode45}} & \multicolumn{3}{c}{\texttt{ode23s}} & \multicolumn{3}{c}{\texttt{ode15s}}\\
\cmidrule(lr){2-4}\cmidrule(lr){5-7}\cmidrule(lr){8-10}
$\delta$ & steps & $|D|$ & error & steps & $|D|$ & error & steps & $|D|$ & error\\
\midrule
$2\cdot10^{-2}$ & 29 & $9.6\cdot10^{-5}$ & $7.0\cdot10^{-5}$ & 64 & $3.2\cdot10^{-2}$ & $4.9\cdot10^{-3}$ & 75 & $6.4\cdot10^{-2}$ & $9.3\cdot10^{-3}$ \\
$10^{-2}$ & 46 & $2.1\cdot10^{-4}$ & $1.1\cdot10^{-4}$ & 72 & $6.3\cdot10^{-2}$ & $9.5\cdot10^{-3}$ & 85 & $1.4\cdot10^{-1}$ & $2.0\cdot10^{-2}$ \\
$10^{-3}$ & 323 & $6.7\cdot10^{-3}$ & $9.8\cdot10^{-4}$ & 89 & $1.4\cdot10^{0}$ & $2.0\cdot10^{-1}$ & 103 & $5.3\cdot10^{0}$ & $6.0\cdot10^{-1}$ \\
$10^{-4}$ & 3040 & $2.2\cdot10^{-1}$ & $3.0\cdot10^{-2}$ & 99 & $2.9\cdot10^{1}$ & $9.5\cdot10^{-1}$ & 140 & $3.3\cdot10^{2}$ & $1.0\cdot10^{0}$ \\
$10^{-5}$ & 30175 & $2.6\cdot10^{0}$ & $3.5\cdot10^{-1}$ & 108 & $2.6\cdot10^{2}$ & $1.0\cdot10^{0}$ & 139 & $1.2\cdot10^{4}$ & $1.0\cdot10^{0}$ \\
$10^{-6}$ & 301475 & $4.3\cdot10^{1}$ & $9.7\cdot10^{-1}$ & 118 & $2.6\cdot10^{3}$ & $1.0\cdot10^{0}$ & 147 & $1.4\cdot10^{5}$ & $1.0\cdot10^{0}$ \\
 
\bottomrule
\end{tabular}
\end{table}

\begin{figure}[t]
\centering
\includegraphics[width=\textwidth]{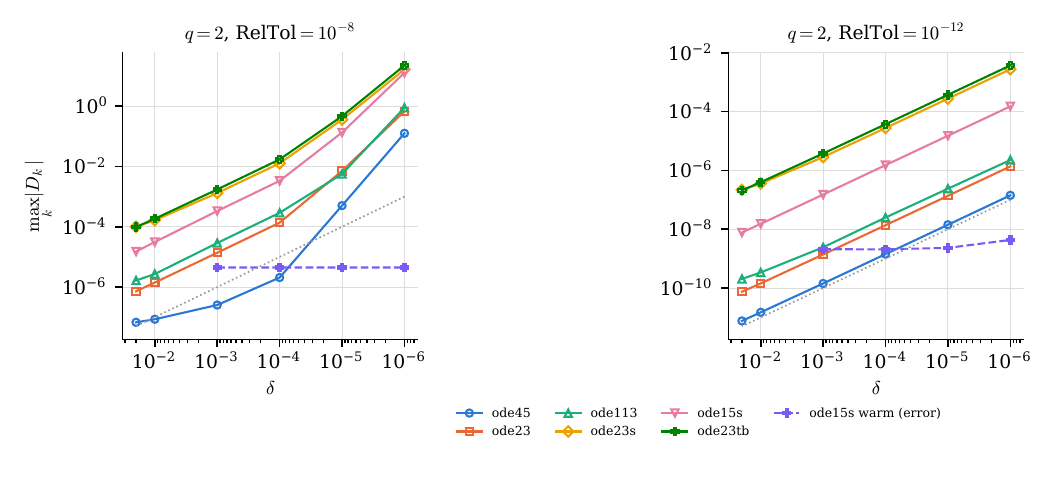}
\caption{Maximal phase delay of the MATLAB solvers versus $\delta$ at tight tolerances, and error of warm-started \texttt{ode15s} (dashed). The dotted line has slope $1/\delta$.}
\label{fig:adaptive}
\end{figure}

\paragraph{Tight tolerances} Table~S1 and Fig.~\ref{fig:adaptive} repeat the experiment with $\mathrm{RelTol}/\mathrm{AbsTol}=10^{-8}/10^{-12}$ and $10^{-12}/10^{-20}$. With $\mathrm{RelTol}=10^{-12}$ (relative regime) the delay of every solver grows by a factor of $10$ per decade of $\delta$, as \eqref{eq:adaptive} predicts. For example, \texttt{ode45} gives $1.4\cdot10^{-10}$, $1.4\cdot10^{-9}$, $1.4\cdot10^{-8}$, $1.4\cdot10^{-7}$ for $\delta=10^{-3},\dots,10^{-6}$. With $\mathrm{RelTol}=10^{-8}$, $\mathrm{AbsTol}=10^{-12}$ the absolute regime takes over once $\delta<10^{-4}$ and the growth steepens towards $\delta^{-(1+\alpha)}$. The constants $K$ differ by more than four orders of magnitude between solvers: the low-order implicit codes \texttt{ode23s}, \texttt{ode23t} and \texttt{ode23tb} have the largest delays. Over all completed runs without warm start in which $\max|D_k|<1$, the ratio of the maximal error to $\max|D_k|$ has median $0.147\approx\frac4{27}$, as Lemma~\ref{lem:phase} predicts when the delay is small compared with the width of the front.

\paragraph{Warm start} Proposition~\ref{prop:translation} gives an exact remedy. Integrate analytically up to a moderate level $y=\eta_0$ (we use $\eta_0=0.01$), i.e.\ start the solver at $t_0=T_\delta(\eta_0)$ with $y(t_0)=\eta_0$. The remaining problem is, up to translation, the same for every $\delta<\eta_0$, and so are the errors. At the benchmark tolerance, warm-started \texttt{ode45}, \texttt{ode23s} and \texttt{ode15s} have maximal errors $1.1\cdot10^{-4}$, $9.5\cdot10^{-3}$ and $2.0\cdot10^{-2}$ for every $\delta\in[10^{-6},10^{-3}]$, with $78$--$83$ and $93$--$103$ steps for the two stiff solvers. At $\mathrm{RelTol}=10^{-8}$ the errors are $1.1\cdot10^{-8}$, $2.0\cdot10^{-5}$ and $4.5\cdot10^{-6}$, again for every $\delta$. The warm start does not change the length of the stiff phase, so explicit codes still need $L_{\mathrm{stiff}}/\beta$ steps. For a general scalar autonomous equation the same is achieved by evaluating $T(\eta_0)$ by quadrature. If the start at $y=\delta$ is kept, \eqref{eq:adaptive} suggests scaling the tolerances with $\delta$ instead; Section~\ref{sec:wp} compares the two remedies.

\subsection{Work and precision}\label{sec:wp}

\begin{figure}[t]
\centering
\includegraphics[width=\textwidth]{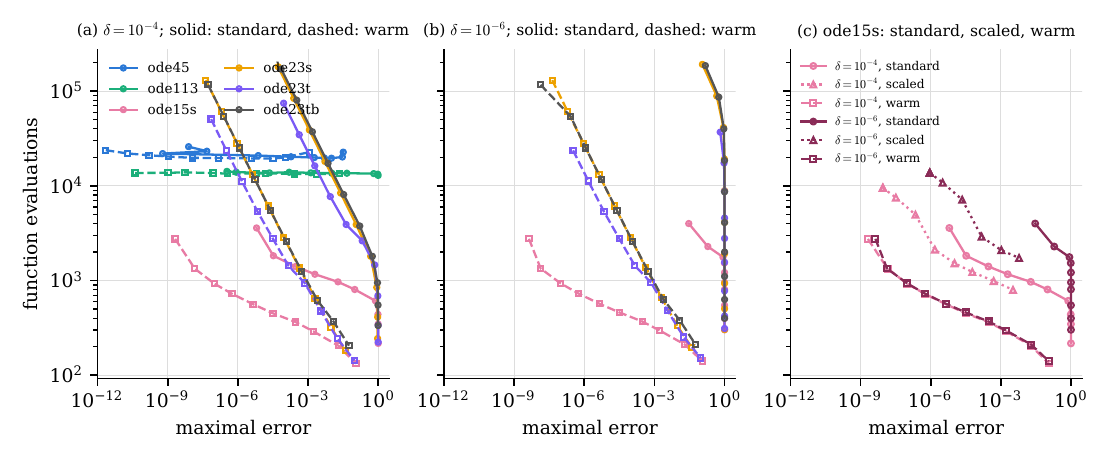}
\caption{Work and precision, MATLAB R2023b, $\mathrm{tol}=10^{-3},\dots,10^{-12}$: function evaluations versus maximal error. (a), (b): standard tolerances (solid) and warm start from $y=0.01$ (dashed); explicit solvers only for $\delta=10^{-4}$. (c): \texttt{ode15s} with standard (solid), scaled (dotted) and warm-started (dashed) integration.}
\label{fig:wp}
\end{figure}

A comparison at fixed tolerance does not say which solver is cheapest for a given accuracy. We therefore measured, with MATLAB R2023b, the number of function evaluations needed to reach a given maximal error, for $\delta=10^{-4}$ and $10^{-6}$ and $\mathrm{tol}=10^{-3},\dots,10^{-12}$, with three strategies:
\emph{standard}, $\mathrm{RelTol}=\mathrm{tol}$, $\mathrm{AbsTol}=\mathrm{tol}/100$ (the ratio of the benchmark setting);
\emph{scaled}, $\mathrm{RelTol}=\mathrm{tol}\cdot\delta$, $\mathrm{AbsTol}=\mathrm{tol}\cdot\delta^2/100$, as suggested by \eqref{eq:adaptive};
and \emph{warm}, the warm start from $y=0.01$ with the standard tolerances (Fig.~\ref{fig:wp}). The explicit solvers were run for $\delta=10^{-4}$ only.
\begin{itemize}
\item With the standard strategy no stiff solver reaches a small error at any tolerance once $\delta$ is small. For $\delta=10^{-4}$ the best errors, at $\mathrm{RelTol}=10^{-12}$, are $6.3\cdot10^{-6}$ for \texttt{ode15s} and $(5$--$9)\cdot10^{-5}$ for \texttt{ode23s}, \texttt{ode23t} and \texttt{ode23tb}; for $\delta=10^{-6}$ they are $0.029$ for \texttt{ode15s} and above $0.11$ for the others.
\item \texttt{ode45} and \texttt{ode113} need $(2.0$--$2.6)\cdot10^4$ and $(1.3$--$1.4)\cdot10^4$ evaluations at every tolerance, since their cost is set by stability in the stiff phase (Section~\ref{sec:adaptive}); their errors still contain the phase delay at loose tolerances but decrease with the tolerance at almost no extra cost.
\item The warm start gives the cheapest solution for every accuracy and removes the dependence on $\delta$. An error of $10^{-6}$ costs \texttt{ode15s} $724$ evaluations for $\delta=10^{-4}$ and $730$ for $\delta=10^{-6}$; an error of $10^{-3}$ costs $364$ (\texttt{ode15s}) to $1354$ (\texttt{ode23s}) evaluations.
\item Scaling the tolerances also makes the error nearly independent of $\delta$ (by \eqref{eq:adaptive} the delay still grows like $\delta^{-(1-\alpha)}$), but at a higher cost. For an error of $10^{-6}$, \texttt{ode15s} needs $4964$ evaluations for $\delta=10^{-4}$ and $13745$ for $\delta=10^{-6}$, $7$ and $19$ times the cost of the warm start. For $\delta=10^{-6}$ the low-order solvers do not reach $10^{-6}$ at all: $\mathrm{tol}\cdot\delta$ falls below the smallest relative tolerance MATLAB accepts ($r_{\min}\approx2.2\cdot10^{-14}$), and \texttt{ode23t} stops near the front with a step-size failure.
\end{itemize}
The SciPy solvers (RK45, DOP853, Radau, BDF, LSODA) give the same picture. With them we also checked the scaled strategy for $\delta=10^{-2},\dots,10^{-6}$ at $\mathrm{tol}=10^{-6}$: the errors of Radau, BDF and LSODA stay within a factor of $5$ over the four decades, while with the standard strategy they grow to $0.09$--$1$. For a scalar problem the warm start is therefore the method of choice. Scaling the tolerances is a valid alternative when $T(\eta_0)$ is not available, at a cost that grows as $\delta$ decreases and with a limit set by the smallest admissible tolerance.

\subsection{How far can the tolerances be tightened?}\label{sec:tol}

\begin{figure}[!tbp]
\centering
\includegraphics[width=\textwidth]{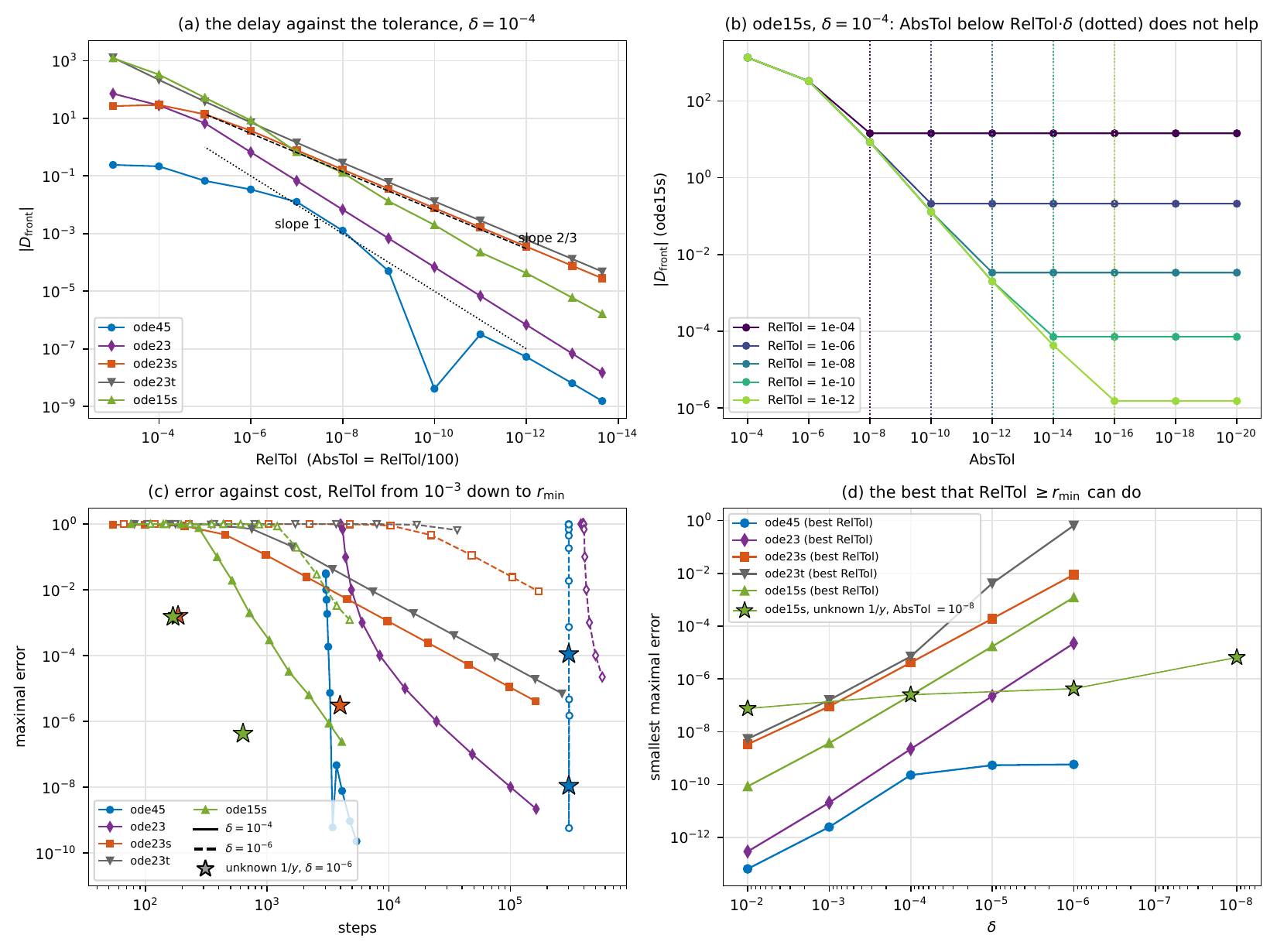}
\caption{Tightening the tolerances, MATLAB R2023b. (a) $|D_{\mathrm{front}}|$ against $\mathrm{RelTol}$, $\mathrm{AbsTol}=\mathrm{RelTol}/100$, $\delta=10^{-4}$. (b) \texttt{ode15s}, $\delta=10^{-4}$: $|D_{\mathrm{front}}|$ against $\mathrm{AbsTol}$ for fixed $\mathrm{RelTol}$; dotted lines at $\mathrm{AbsTol}=\mathrm{RelTol}\cdot\delta$. (c) Maximal error against the number of accepted steps, $\mathrm{RelTol}=10^{-3},\dots,r_{\min}$, $\delta=10^{-4}$ (solid) and $10^{-6}$ (dashed); stars: unknown $1/y$ \eqref{eq:clock}, $\mathrm{RelTol}=r_{\min}$, $\mathrm{AbsTol}=10^{-4}$ and $10^{-8}$, $\delta=10^{-6}$. (d) Smallest maximal error reached with $\mathrm{RelTol}\ge r_{\min}$ against $\delta$; stars: \texttt{ode15s} with the unknown $1/y$ and $\mathrm{AbsTol}=10^{-8}$.}
\label{fig:tol}
\end{figure}

By \eqref{eq:adaptive} the delay can be made small by decreasing the tolerances. To see at what price, we ran \texttt{ode45}, \texttt{ode23}, \texttt{ode23s}, \texttt{ode23t} and \texttt{ode15s} for $\delta=10^{-2},\dots,10^{-6}$ with $\mathrm{RelTol}=10^{-3},\dots,10^{-13}$ and $r_{\min}$, always with $\mathrm{AbsTol}=\mathrm{RelTol}/100$, and for $\delta=10^{-4}$ with the two tolerances varied independently (Fig.~\ref{fig:tol}, Table~\ref{tab:tolneeded}).

\paragraph{The delay decreases like a power of the tolerance} For $\delta=10^{-4}$ and $10^{-11}\le\mathrm{RelTol}\le10^{-6}$ the slopes of $|D_{\mathrm{front}}|$ are $0.67$ (\texttt{ode23s}), $0.68$ (\texttt{ode23t}), $0.90$ (\texttt{ode15s}) and $1.00$ (\texttt{ode23}), as \eqref{eq:adaptive} predicts (Fig.~\ref{fig:tol}a). The small delay of \texttt{ode45} changes sign, and below $\mathrm{RelTol}\approx10^{-10}$ it is set by rounding.

\paragraph{The delay grows faster than $1/\delta$} At $\mathrm{RelTol}=10^{-8}$ the delay grows like $\delta^{-1.67}$ (\texttt{ode23s}) to $\delta^{-1.99}$ (\texttt{ode23}) for $10^{-6}\le\delta\le10^{-3}$, because with $\mathrm{AbsTol}=\mathrm{RelTol}/100$ the absolute test is active at the start (second case of \eqref{eq:adaptive}). To keep the delay fixed, the tolerance must decrease by about $2.5$ decades per decade of $\delta$ for the second-order stiff codes, and by about $2$ decades for \texttt{ode15s} and \texttt{ode23}.

\paragraph{Both tolerances must go down} At $y=\delta$ the error test admits $\max(\mathrm{RelTol}\cdot\delta,\mathrm{AbsTol})$. Decreasing $\mathrm{AbsTol}$ reduces the delay only until $\mathrm{AbsTol}=\mathrm{RelTol}\cdot\delta$ (Fig.~\ref{fig:tol}b). With $\mathrm{RelTol}=10^{-4}$ and $\mathrm{AbsTol}=10^{-20}$ the fronts of \texttt{ode15s} and \texttt{ode23s} still come $14.3$ and $6.1$ time units early.

\paragraph{The price and the limit} Table~\ref{tab:tolneeded} gives, for a maximal error of at most $10^{-2}$, the largest $\mathrm{RelTol}$ that achieves it. For the stiff codes the cost rises quickly: at $\delta=10^{-6}$ \texttt{ode15s} needs $25$ times its benchmark number of steps, \texttt{ode23s} $1446$ times, and \texttt{ode23t} does not reach an error of $10^{-2}$ at all; it stops near the front with a step-size failure for $\delta\le10^{-5}$ and $\mathrm{RelTol}\le10^{-13}$. A maximal error of $10^{-4}$ is out of reach for all three stiff codes at $\delta=10^{-6}$, and for \texttt{ode23s} and \texttt{ode23t} already at $\delta=10^{-5}$. The smallest maximal error that any admissible tolerance gives grows steadily as $\delta$ decreases (Fig.~\ref{fig:tol}d); at $\delta=10^{-6}$ it is $1.2\cdot10^{-3}$ for \texttt{ode15s}, $8.9\cdot10^{-3}$ for \texttt{ode23s} and $0.65$ for \texttt{ode23t}. For the explicit codes tightening is almost free, since their cost is set by the stiff phase (\texttt{ode45} reaches an error of $10^{-4}$ at $\delta=10^{-4}$ with $3282$ steps instead of $3040$), but this cost grows like $1/\delta$. In comparison, the clock variable of Section~\ref{sec:sundman} gives \texttt{ode15s} at $\delta=10^{-6}$ a maximal error of $1.5\cdot10^{-3}$ with $168$ steps and of $4.3\cdot10^{-7}$ with $633$ steps (stars in Fig.~\ref{fig:tol}c,\,d). Tightening the tolerances thus works in principle, but for the stiff codes it is practical only down to about $\delta=10^{-4}$.

\begin{table}[!tbp]
\centering\footnotesize
\setlength{\tabcolsep}{3.5pt}
\caption{The price of tightening the tolerances (MATLAB R2023b, $\mathrm{AbsTol}=\mathrm{RelTol}/100$): the largest $\mathrm{RelTol}$ that gives a maximal error of at most $10^{-2}$, the number of accepted steps and, in brackets, the number of steps relative to the benchmark tolerance.}
\label{tab:tolneeded}
\begin{tabular}{lccccc}
\toprule
$\delta$ & \texttt{ode45} & \texttt{ode23} & \texttt{ode23s} & \texttt{ode23t} & \texttt{ode15s}\\
\midrule
$10^{-2}$ & $10^{-3}$: 39 (0.8) & $10^{-4}$: 95 (1.0) & $10^{-4}$: 72 (1.0) & $10^{-5}$: 241 (2.1) & $10^{-5}$: 133 (1.6)\\
$10^{-3}$ & $10^{-3}$: 317 (1.0) & $10^{-5}$: 572 (1.2) & $10^{-7}$: 891 (10) & $10^{-7}$: 1433 (9.9) & $10^{-7}$: 318 (3.1)\\
$10^{-4}$ & $10^{-6}$: 3074 (1.0) & $10^{-8}$: 6033 (1.5) & $10^{-9}$: 4525 (46) & $10^{-9}$: 7411 (47) & $10^{-9}$: 716 (5.1)\\
$10^{-5}$ & $10^{-9}$: 30418 (1.0) & $10^{-10}$: 49851 (1.3) & $10^{-12}$: 47101 (436) & $10^{-12}$: 78089 (462) & $10^{-11}$: 1649 (12)\\
$10^{-6}$ & $10^{-11}$: 302087 (1.0) & $10^{-12}$: 445711 (1.1) & $r_{\min}$: 170639 (1446) & none & $10^{-13}$: 3729 (25)\\
 
\bottomrule
\end{tabular}
\end{table}

\subsection{A remedy without the exact solution: time as a dependent variable}\label{sec:sundman}

The warm start uses $T(\eta_0)$, and scaling or tightening the tolerances is expensive. The analysis above suggests a third remedy that needs neither. The delay arises because the codes control the error in $y$, whereas the quantity that is sensitive in the slow phase is the time at which $y$ is reached. We therefore make time a dependent variable. With a positive rate function $\rho(y)$ and a new independent variable $\tau$ defined by $\dd t/\dd\tau=1/\rho(y)$ (a Sundman transformation \cite{Sundman1912,Hairer1997}), \eqref{eq:fq} becomes
\begin{equation}\label{eq:sundman}
\frac{\dd y}{\dd\tau}=\frac{f_q(y)}{\rho(y)},\qquad \frac{\dd t}{\dd\tau}=\frac{1}{\rho(y)},\qquad y(0)=\delta,\quad t(0)=0 .
\end{equation}
For $\rho(y)=y^{q-1}$, the scale-invariant choice for $y'\approx y^q$ \cite{BuddPiggott2001}, the first equation is the logistic equation $\dd y/\dd\tau=y(1-y)$: the degeneracy at $y=0$ and with it the long slow phase disappear, and the front is reached after $\tau\approx\ln(1/\delta)$. The rate function $\rho(y)=|f_q(y)|/|y|+|f_q'(y)|$, which is positive on $(0,1]$, achieves the same without knowledge of $q$.

The essential ingredient is the error control: we impose a relative tolerance on $y$ and an \emph{absolute} tolerance on $t$. The local error of the $t$-component over a step of length $k$ in $\tau$ is of order $k^{p+1}y^{1-q}$, so the absolute control forces $k^{p+1}\lesssim\mathrm{tol}\,y^{q-1}$; then the local error in $y$, of order $k^{p+1}y$, is at most $\mathrm{tol}\,y^q$, and the corresponding time defect $e_y/f_q(y)$ is of order $\mathrm{tol}$ in every step. The step size is thus reduced exactly where the timing is sensitive, by an amount that grows only like $\delta^{-(q-1)/(p+1)}$, and only for the steps near $y=\delta$. The defects still add up over the slow phase, whose number of steps grows like $(\mathrm{tol}\,\delta^{q-1})^{-1/(p+1)}$, so the delay behaves like $\mathrm{tol}^{\alpha}\delta^{-(q-1)(1-\alpha)}$ with $\alpha=p/(p+1)$, instead of $\mathrm{tol}^{\alpha}\delta^{1-q}$. For the adaptive Euler method with $\rho=y$ and the local error of $t$ equal to $\mathrm{tol}$, each step adds about $2\,\mathrm{tol}$ to the delay and there are about $\sqrt{2/(\mathrm{tol}\,\delta)}$ steps, so $D\approx\sqrt{8\,\mathrm{tol}/\delta}$ instead of $\sqrt{\mathrm{tol}}/\delta$ (Fig.~\ref{fig:euler}b); the computed delays are within $1.5\%$ of this value for $\delta\le10^{-4}$, and the delay is now spread over the whole slow phase (Fig.~\ref{fig:euler}a). For a first-order method the gain is thus only a factor $\delta^{-1/2}$, at the price of many steps; for codes of higher order the remaining growth is weak. A blow-up time has been estimated with a Sundman transformation in \cite{StuartFloater1990}, and step sizes have recently been chosen from the sensitivity of a hitting time \cite{HoelMeo2026}. Here the absolute control of $t$ substantially reduces the sensitivity of the front timing without requiring the exact solution, although the model above still predicts a residual dependence on $\delta$, especially for low-order methods.

Rescalings of this kind are standard for blow-up problems, where they go back to the rescaling algorithm of Berger and Kohn \cite{BergerKohn1988}. A rescaling of the variables alone does not help, however. The linear change $z=y/\delta$, $s=\delta t$ of Section~\ref{sec:structure} turns the slow phase into an interval of length $O(1)$, but a Runge--Kutta code with relative error control produces essentially the same points in $(s,z)$ as in $(t,y)$, and the front is displaced by the same amount: with RK45 at $\mathrm{tol}=10^{-6}$ the delays of the two formulations agree to within $0.5\%$ for $\delta\le10^{-4}$. What reduces the delay in \eqref{eq:sundman} is the absolute control of $t$. Combustion codes reduce stiffness by a different mechanism, removing fast modes with quasi-steady-state and partial-equilibrium approximations \cite{MottOranvanLeer2000,LuLawYooChen2009}; this lowers the cost of fast modes in large reaction mechanisms but does not address the timing of ignition.

Table~S2 compares the standard formulation with \eqref{eq:sundman}, both solved with LSODA \cite{Petzold1983} as implemented in SciPy, with $\mathrm{RelTol}=\mathrm{tol}$ on $y$ and $\mathrm{AbsTol}=\mathrm{tol}$ on $t$. With the standard formulation the error is close to $1$ for $\delta\le10^{-6}$ ($q=2$) and $\delta\le10^{-4}$ ($q=3$) even at $\mathrm{tol}=10^{-8}$. With time as a dependent variable the front delay is of the order of $\mathrm{tol}$, except in the two cases discussed below, and the error is at most $1.1\cdot10^{-4}$ at $\mathrm{tol}=10^{-4}$ and at most $2.2\cdot10^{-6}$ at $\mathrm{tol}=10^{-8}$ for every $\delta$, at a cost that is lower than that of the standard formulation in all but one case. The generic rate function costs up to twice as much, with errors within a factor of $20$ of those for $\rho=y^{q-1}$. The largest errors at $\mathrm{tol}=10^{-8}$ occur for $\delta=10^{-8}$ ($q=2$) and $\delta=10^{-4}$ ($q=3$), where $t\approx10^{8}$ near the front. There the relative part of the test on $t$, $10^{-13}t\approx10^{-5}$, exceeds $\mathrm{tol}$, and $t$ itself is represented only to about $u\,t\approx10^{-8}$: these two runs are limited by both. If $t$ is given only the relative tolerance, as in codes with a scalar $\mathrm{RelTol}$, the error is still reduced by one to two orders of magnitude but grows again as $\delta$ decreases, since $t$ itself is of order $\delta^{1-q}$; this is the situation of the MATLAB solvers (see below). In SUNDIALS \cite{Hindmarsh2005} the absolute control of $t$ is available without modifying the code, through a user-supplied error weight function (\texttt{CVodeWFtolerances}). With CVODE (BDF) and the weight $1/(10^{-13}|t|+\mathrm{tol})$ for $t$, the setting used with LSODA, the maximal error is at most $5.5\cdot10^{-3}$ at $\mathrm{tol}=10^{-4}$ and at most $3.5\cdot10^{-6}$ at $\mathrm{tol}=10^{-8}$ apart from the same two cases, two to five orders of magnitude below that of the standard formulation, at $1.0$--$2.0$ times its cost (except for one run in which the standard formulation skips the front in $10$ function evaluations). The standard tolerance interface, with one scalar relative tolerance for both components, leaves errors above $0.75$ at $\mathrm{tol}=10^{-4}$ for every $\delta\le10^{-4}$ ($q=2$) and every $\delta$ ($q=3$). Unlike with LSODA, the CVODE delay still grows as $\delta$ decreases, as the accumulation of the per-step defects described above predicts.

\paragraph{MATLAB} The error test of the MATLAB solvers, $|e_i|\le\max(\mathrm{RelTol}\,|u_i|,\mathrm{AbsTol}_i)$, uses one scalar relative tolerance for all components $u_i$, so for \eqref{eq:sundman} in the unknowns $(y,t)$ the tolerance on $t$ is relative. A change of variable removes this restriction. With $z=\ln y$, \eqref{eq:sundman} becomes
\begin{equation}\label{eq:sundmanlog}
\frac{\dd z}{\dd\tau}=\frac{f_q(y)}{y\,\rho(y)},\qquad \frac{\dd t}{\dd\tau}=\frac{1}{\rho(y)},\qquad y=e^{z},
\end{equation}
and an absolute tolerance on $z$ is a relative tolerance on $y$. With $\mathrm{RelTol}$ at its smallest admissible value $r_{\min}\approx2.2\cdot10^{-14}$ and $\mathrm{AbsTol}=\mathrm{tol}$ for both components, the unmodified solvers therefore control $y$ relatively and $t$ absolutely, as long as $r_{\min}|z|$ and $r_{\min}|t|$ stay below $\mathrm{tol}$; the exceptions are discussed below. Table~S3 compares the three formulations for \texttt{ode15s} and \texttt{ode23s} (MATLAB R2023b, the problems and tolerances of Table~S2). In the variables $(y,t)$ the errors of both codes are still above $0.85$ at $\mathrm{tol}=10^{-4}$ for $\delta\le10^{-6}$ ($q=2$) and $\delta\le10^{-3}$ ($q=3$). With \eqref{eq:sundmanlog} the maximal error of \texttt{ode15s} is at most $2.2\cdot10^{-3}$ at $\mathrm{tol}=10^{-4}$ and at most $8.0\cdot10^{-7}$ at $\mathrm{tol}=10^{-8}$ apart from the two cases with $t\approx10^8$, where $r_{\min}t\approx2\cdot10^{-6}$ exceeds $\mathrm{tol}$, at $0.75$--$2.0$ times the cost of the standard formulation, whose errors are close to $1$ for the smaller values of $\delta$. The second-order codes \texttt{ode23s}, \texttt{ode23t} and \texttt{ode23tb} also improve, to errors between $1.2\cdot10^{-6}$ and $9.8\cdot10^{-2}$, but at $1.5$--$117$ times the cost, and their delays grow as $\delta$ decreases. This is the accumulation of per-step defects described above: for $\alpha=2/3$ it predicts $D\propto\mathrm{tol}^{2/3}\delta^{-(q-1)/3}$, and the measured exponents of $|D|$ are $0.62$--$0.67$ in $\mathrm{tol}$ and $-0.34$ to $-0.40$ ($q=2$) and $-0.67$ to $-0.72$ ($q=3$) in $\delta$. For \texttt{ode15s}, whose order varies up to $5$, they are $0.75$--$0.80$ in $\mathrm{tol}$ and $-0.27$ to $-0.31$ ($q=2$) and $-0.39$ to $-0.46$ ($q=3$) in $\delta$, between the values for $\alpha=2/3$ and $\alpha=5/6$. The explicit codes, whose cost is set by the stiff phase, reach errors of at most $1.1\cdot10^{-4}$ (\texttt{ode45}) and $1.8\cdot10^{-3}$ (\texttt{ode113}) at $\mathrm{tol}=10^{-4}$, at $0.74$--$1.07$ times the cost of the standard formulation; three runs of \texttt{ode45} at $\mathrm{tol}=10^{-4}$ stopped after $10$--$12$ steps without reaching $t=2/((q-1)\delta^{q-1})$ and are not included. The generic rate function gives errors between $0.14$ and $1.03$ times those for $\rho=y^{q-1}$, at $0.95$--$1.84$ times the cost.

\paragraph{A clock as unknown} The time can also be controlled without changing the independent variable. With $w=y^{1-q}$, \eqref{eq:fq} becomes
\begin{equation}\label{eq:clock}
w'=-(q-1)\bigl(1-w^{-1/(q-1)}\bigr),\qquad w(0)=\delta^{1-q},
\end{equation}
and in the slow phase $w\approx\delta^{1-q}-(q-1)t$: the unknown is a clock, and an absolute error $\varepsilon$ in $w$ is a time shift $\varepsilon/(q-1)$ in the slow phase (in general $\varepsilon/((q-1)(1-y))$, since $w'=-(q-1)(1-y)$). An absolute tolerance on $w$ is therefore an absolute tolerance on the time at which each level is reached, whereas a relative tolerance on $w$ is a relative tolerance on $y$ and gains nothing. This is also the variable in which explicit Euler has a $\delta$-independent delay (Remark~\ref{rem:variables}). With $\mathrm{RelTol}=r_{\min}$ and $\mathrm{AbsTol}=\mathrm{tol}$ the unmodified MATLAB solvers integrate \eqref{eq:clock} in the original time variable; output at prescribed times is then available without event location. Table~\ref{tab:clock} shows the result for the problems and tolerances of Table~S3. Over all six solvers the maximal error is at most $2.7\cdot10^{-3}$ at $\mathrm{tol}=10^{-4}$ and $1.7\cdot10^{-5}$ at $\mathrm{tol}=10^{-8}$ for $q=2$, down to $\delta=10^{-8}$, and at most $3.7\cdot10^{-3}$ and $6.0\cdot10^{-5}$ for $q=3$; the largest values occur for $\mathrm{tol}=10^{-8}$ and $w(0)=10^{8}$ ($q=2$, $\delta=10^{-8}$ and $q=3$, $\delta=10^{-4}$). There $r_{\min}w(0)\approx2.2\cdot10^{-6}$ exceeds $\mathrm{tol}$, so the test is relative at the start, and $\mathrm{tol}$ is also close to the rounding level $u\,w(0)\approx1.1\cdot10^{-8}$; these runs are limited by both. The median cost equals that of the standard formulation (ratio $1.03$, range $0.72$--$14$, the largest for \texttt{ode15s} at $\delta=10^{-8}$, $\mathrm{tol}=10^{-4}$). In $60$ of the $73$ runs that can be compared, \eqref{eq:clock} is more accurate than \eqref{eq:sundmanlog}, mainly for the second-order codes: \texttt{ode23s} with $q=2$, $\delta=10^{-8}$, $\mathrm{tol}=10^{-4}$ has error $2.1\cdot10^{-3}$ at $1.8$ times the cost of the standard formulation, against $9.8\cdot10^{-2}$ at $120$ times with \eqref{eq:sundmanlog}. With $\mathrm{RelTol}=\mathrm{AbsTol}=\mathrm{tol}$ instead, the errors for \eqref{eq:clock} reach $0.98$ ($q=3$) and $1.0$ ($q=2$), and two runs of \texttt{ode23s} overflow.

\begin{table}[t]
\centering\small
\setlength{\tabcolsep}{4pt}
\caption{MATLAB R2023b with the unknown $w=y^{1-q}$ \eqref{eq:clock}, $\mathrm{RelTol}=r_{\min}$, $\mathrm{AbsTol}=\mathrm{tol}$, for the problems of Table~S3: maximal error at the computed points and cost (function evaluations divided by those of the standard formulation). Explicit solvers only for $t_{\mathrm{end}}\le2\cdot10^{6}$.}
\label{tab:clock}
\begin{tabular}{ccc cc cc cc}
\toprule
 & & & \multicolumn{2}{c}{\texttt{ode45}} & \multicolumn{2}{c}{\texttt{ode15s}} & \multicolumn{2}{c}{\texttt{ode23s}}\\
\cmidrule(lr){4-5}\cmidrule(lr){6-7}\cmidrule(lr){8-9}
$q$ & $\delta$ & tol & error & cost & error & cost & error & cost\\
\midrule
$2$ & $10^{-2}$ & $10^{-4}$ & $1.1\cdot10^{-4}$ & 0.9 & $4.7\cdot10^{-4}$ & 0.9 & $6.6\cdot10^{-4}$ & 1.1 \\
 & $10^{-2}$ & $10^{-8}$ & $1.1\cdot10^{-8}$ & 0.8 & $7.6\cdot10^{-8}$ & 0.9 & $1.3\cdot10^{-6}$ & 1.1 \\
 & $10^{-4}$ & $10^{-4}$ & $1.1\cdot10^{-4}$ & 1.0 & $1.0\cdot10^{-3}$ & 0.7 & $1.1\cdot10^{-3}$ & 1.3 \\
 & $10^{-4}$ & $10^{-8}$ & $1.1\cdot10^{-8}$ & 1.0 & $2.5\cdot10^{-7}$ & 0.9 & $2.1\cdot10^{-6}$ & 1.3 \\
 & $10^{-6}$ & $10^{-4}$ & $1.1\cdot10^{-4}$ & 1.0 & $1.5\cdot10^{-3}$ & 0.9 & $1.6\cdot10^{-3}$ & 1.5 \\
 & $10^{-6}$ & $10^{-8}$ & $1.1\cdot10^{-8}$ & 1.0 & $4.3\cdot10^{-7}$ & 1.0 & $3.1\cdot10^{-6}$ & 1.8 \\
 & $10^{-8}$ & $10^{-4}$ & -- & -- & $2.1\cdot10^{-3}$ & 14.4 & $2.1\cdot10^{-3}$ & 1.8 \\
 & $10^{-8}$ & $10^{-8}$ & -- & -- & $6.5\cdot10^{-6}$ & 1.1 & $1.4\cdot10^{-5}$ & 2.0 \\
\midrule
$3$ & $10^{-2}$ & $10^{-4}$ & $5.6\cdot10^{-5}$ & 1.0 & $3.2\cdot10^{-4}$ & 1.1 & $7.0\cdot10^{-4}$ & 2.9 \\
 & $10^{-2}$ & $10^{-8}$ & $5.6\cdot10^{-9}$ & 1.0 & $1.1\cdot10^{-7}$ & 1.2 & $1.4\cdot10^{-6}$ & 2.9 \\
 & $10^{-3}$ & $10^{-4}$ & $5.5\cdot10^{-5}$ & 1.0 & $6.2\cdot10^{-4}$ & 0.9 & $1.5\cdot10^{-3}$ & 4.8 \\
 & $10^{-3}$ & $10^{-8}$ & $5.3\cdot10^{-9}$ & 1.0 & $2.4\cdot10^{-7}$ & 1.5 & $3.1\cdot10^{-6}$ & 5.0 \\
 & $10^{-4}$ & $10^{-4}$ & -- & -- & $7.8\cdot10^{-4}$ & 1.1 & $3.1\cdot10^{-3}$ & 8.8 \\
 & $10^{-4}$ & $10^{-8}$ & -- & -- & $5.3\cdot10^{-6}$ & 1.6 & $4.9\cdot10^{-5}$ & 6.5 \\
 
\bottomrule
\end{tabular}
\end{table}

\paragraph{Nodal and common-grid errors} Except in this paragraph, errors are maximal nodal errors, taken over the accepted steps, which differ between methods and formulations. To check that the comparisons do not depend on this sampling, we evaluated the standard and the clock formulation on a common set of output times with the solvers' own interpolants (MATLAB \texttt{deval}): $2001$ uniform points on $[0,t_{\mathrm{end}}]$ and $4001$ points on each of the windows of length $80$ around the exact and the computed front. The set does not contain the accepted steps; for the clock, the interpolant of $w$ is evaluated and inverted. The $52$ runs, all completed, are $q=2$, $\delta=10^{-4},10^{-6}$ and $q=3$, $\delta=10^{-3},10^{-4}$, with $\mathrm{tol}=10^{-4},10^{-8}$, both formulations and \texttt{ode23s}, \texttt{ode23t}, \texttt{ode15s}, and \texttt{ode45} only for $q=2$, $\delta=10^{-4}$ (Table~S6). With four times as many output times the maxima change by at most 0.02\%. For the standard formulation the common-grid error exceeds the nodal error by at most $8\%$. For the clock the same holds in $25$ of the $26$ runs; the exception is \texttt{ode23t} for $q=3$, $\delta=10^{-4}$, $\mathrm{tol}=10^{-4}$, whose interpolant gives an error of $2.5\cdot10^{-2}$ between accepted steps, against a nodal error of $3.7\cdot10^{-3}$ and an error of $1.0$ for the standard formulation. Continuous output of the clock formulation is therefore not uniformly as accurate as its nodal values. The same check with the dense output of SciPy's BDF, Radau, LSODA and RK45 gives ratios of at most $1.09$ (standard) and $1.54$ (clock), and the maxima change by at most 1.6\% under the same refinement. The differences of several orders of magnitude between the formulations are therefore not an artefact of the sampling.

\paragraph{Other implementations} We repeated the experiment for $q=2$ with the solvers of the Julia package OrdinaryDiffEq \cite{RackauckasNie2017} (Julia 1.13.1, OrdinaryDiffEq 7.8.1, OrdinaryDiffEqBDF 2.4.12, OrdinaryDiffEqSDIRK 2.9.7), most of whose methods accept component-wise relative and absolute tolerances (exceptions below). Both formulations solve the same initial value problem, $y'=y^2-y^3$, $y(0)=\delta$, $0\le t\le2/\delta$; the errors are measured against its exact solution at the computed points $(t_k,y_k)$. We used $\mathrm{RelTol}=10^{-6}$, $\mathrm{AbsTol}=10^{-8}$ for $y$ throughout and, in the transformed formulation, $\mathrm{AbsTol}_t=10^{-4},10^{-6},10^{-8}$ with $\mathrm{RelTol}_t=10^{-13}$ for $t$ (Table~S4). For the explicit methods Tsit5 and Vern7 the delay follows $\mathrm{AbsTol}_t$ closely at almost no extra cost, since their cost is set by the stiff phase. For QNDF, a variable-order NDF method like \texttt{ode15s}, the front at $\delta=10^{-6}$ is displaced by $2.5\cdot10^{4}$ time units in the standard formulation and by $1.9\cdot10^{-4}$ with $\mathrm{AbsTol}_t=10^{-6}$, at $2428$ instead of $962$ function evaluations. The low-order implicit methods TRBDF2 and KenCarp4 respond as well, with a floor of $10^{-4}$--$10^{-3}$ at $\delta=10^{-6}$. The remedy requires that the solver honors the tolerance of each component in its error estimate. FBDF does not respond to $\mathrm{AbsTol}_t$, and the Rosenbrock methods (Rosenbrock23, the method of \texttt{ode23s}, and Rodas5P) do not accept a vector of relative tolerances in these versions.

\paragraph{Which problems need it} The delay law and the remedies use two properties of the family \eqref{eq:fq}: the solution leaves a degenerate equilibrium, $f(y)\approx c\,y^q$ with $q\ge2$, and the time of a later front is the quantity of interest. For a scalar autonomous equation $x'=f(x)$ with these properties the same computations apply with the local time scale
\[
\theta(x)=\frac{|x|}{|f(x)|}\approx\frac{x^{1-q}}{c}
\]
in place of $y^{1-q}$: by the integral of Proposition~\ref{prop:defects}, relative error control keeps the time defect of a step of the order $\mathrm{RelTol}\cdot\theta(x)$, which is unbounded near the equilibrium, and the defects add up over the steps, as in \eqref{eq:adaptive}. Time as a dependent variable with an absolute tolerance on $t$ and the generic rate function $\rho=1/\theta+|f'|$, which near a degenerate equilibrium equals $(q+1)/\theta$ to leading order, or the clock $x^{1-q}$ with an absolute tolerance, then control the time. Other ignition and autocatalytic models and blow-up problems, where the blow-up time is the quantity of interest \cite{StuartFloater1990,HoelMeo2026}, are of this kind. The mechanism itself needs only that the equation is autonomous. Other ways in which a slow passage arises (exponential growth from data below the absolute tolerance, a bottleneck, recurrent passages), systems, and a condition number that tells in advance whether a given error test protects the time of an event are treated in \cite{BorosC2026}.

\subsection{The remedy in practice}\label{sec:recipe}

For a scalar equation $y'=f(y)$ whose solution leaves a degenerate equilibrium $y=0$, $f(y)\approx c\,y^q$ with $q\ge2$, as in \eqref{eq:fq}, the clock remedy of Section~\ref{sec:sundman} reads as follows.
\begin{enumerate}
\item \emph{Unknown.} Integrate for $w=y^{1-q}$, i.e.\ $w'=(1-q)\,y^{-q}f(y)$ with $y=w^{-1/(q-1)}$; for \eqref{eq:fq} this is \eqref{eq:clock}, and for the flame problem $w=1/y$, $w'=-(1-1/w)$. Near $y=0$, $w'\approx-(q-1)c$: the unknown is a clock with speed $(q-1)c$. The Jacobian is $\partial w'/\partial w=-w^{-q/(q-1)}$ for \eqref{eq:clock} ($-1/w^2$ for $q=2$) and should be supplied to implicit codes. If $q$ is not known, use time as a dependent variable \eqref{eq:sundmanlog} with the generic rate function instead.
\item \emph{Interval and output.} Keep the original independent variable, the initial value $w(0)=\delta^{1-q}$ and the original interval; recover $y=w^{-1/(q-1)}$ at the output points. No event location is needed.
\item \emph{Tolerances.} $\mathrm{AbsTol}=\mathrm{tol}$ on $w$, which in the slow phase controls the time at which each level is reached to about $\mathrm{tol}/((q-1)|c|)$ per step ($\mathrm{tol}/(q-1)$ for \eqref{eq:fq}; a per-step interpretation, not a bound on the accumulated error), and the smallest admissible relative tolerance ($r_{\min}$ in MATLAB), so that the error test is absolute. The relative part of the test, $r_{\min}|w|\le r_{\min}\delta^{1-q}$, must stay below $\mathrm{tol}$, i.e.\ $\mathrm{tol}\gtrsim r_{\min}\delta^{1-q}$; otherwise the test is relative at the start, as in the runs with $w(0)=10^8$ and $\mathrm{tol}=10^{-8}$ of Table~\ref{tab:clock}. This bound is about $200$ times the rounding level $u\,\delta^{1-q}$ of $w(0)$ and of the front time $t\approx\delta^{1-q}/(q-1)$. In codes with componentwise relative tolerances the same choice applies to $w$ alone.
\item \emph{Safeguards.} For \eqref{eq:fq} the exact $w$ decreases monotonically to $1$ and stays positive; after the front $w'\approx-(w-1)$, the same stiffness as for $y$. A computed $w\le0$ lies outside the positive branch on which $y=w^{-1/(q-1)}$ represents the solution (for some exponents an algebraic inverse exists, but it is not the solution studied); it should make the code reject the step or stop with a message; it can result from a step beyond the stability boundary, from a failed nonlinear solve or from extrapolation in a predictor. Our scripts evaluate the right-hand side with $\max(w,\mathrm{realmin})$; all $76$ MATLAB runs with this unknown and the tolerances of item~3 completed. With $\mathrm{RelTol}=\mathrm{AbsTol}$, in contrast, the relative test dominates and the remedy fails (Section~\ref{sec:sundman}).
\item \emph{Cost.} In the median about that of the standard formulation (ratio $1.03$, range $0.72$--$14$ over the runs of Table~\ref{tab:clock}).
\end{enumerate}
The three remedies differ in what can be claimed. The warm start is exact: the remaining problem is, up to translation, independent of $\delta$, apart from rounding. Time as a dependent variable reduces the sensitivity of the timing, but the model of Section~\ref{sec:sundman} and the experiments retain a residual dependence $D\propto\mathrm{tol}^{\alpha}\delta^{-(q-1)(1-\alpha)}$ for the low-order codes. The clock unknown gave errors that are nearly independent of $\delta$ over the range tested, within the limits set by $r_{\min}$ and by rounding; we have no global error theorem for it. The recipe needs a scalar unknown that does not change sign. For systems, the corresponding unknowns are those that run like clocks during the slow motion, such as logarithms of small positive components \cite{BorosC2026}; a general construction is beyond the scope of this paper.

%=====================================================================
\section{A second example: $y'=y^3(1-y)$}\label{sec:q3}

In chemical kinetics, $f_q(y)=y^q(1-y)$ is the rate of the autocatalytic step $\mathrm{A}+q\,\mathrm{B}\to(q+1)\,\mathrm{B}$ in a closed system with normalized total concentration \cite{GrayScott1990}. The case $q=2$ (cubic autocatalysis) coincides with the flame model; we now take $q=3$, which has a stronger degeneracy at $y=0$:
\begin{equation}\label{eq:q3}
y'=y^3(1-y),\qquad y(0)=\delta,\qquad 0\le t\le 1/\delta^2 .
\end{equation}
The slow phase now lasts about $1/(2\delta^2)$; for $\delta=10^{-4}$ this is $5\cdot10^7$ time units.

\subsection{Structure}
Every statement of Section~\ref{sec:structure} has an analogue for \eqref{eq:fq}, with the same proofs.
\begin{itemize}
\item \emph{Translation.} $y_\delta(t)=Y_q\bigl(t+\Theta_q(\delta)\bigr)$ with the master profile $Y_q:=\Theta_q^{-1}$ from \eqref{eq:Theta}. For $q=3$, $Y_3$ is not elementary; we evaluate it by bisection in the logit variable $\ln(y/(1-y))$ followed by Newton's method.
\item \emph{Sensitivity.} A relative perturbation $\rho$ of the initial value shifts the solution by $\rho\,\delta/f_q(\delta)\approx\rho/\delta^{q-1}$.
\end{itemize}

\paragraph{Other nonlinearities} The same holds for $y'=y^q(1-y^\ell)$ with integers $q\ge2$, $\ell\ge1$. Translation invariance and Lemma~\ref{lem:phase} use only that the equation is scalar and autonomous, and the scaling argument behind Theorem~\ref{prop:delay} depends only on the degeneracy $y^q$ at $y=0$, so the heuristic of Section~\ref{sec:adaptive} again suggests $D\approx K\,\mathrm{RelTol}^{\alpha}/\delta^{q-1}$ in the relative regime. What changes is the time function. Since $1/\bigl(y^q(1-y^\ell)\bigr)=\sum_{j\ge0}y^{j\ell-q}$ for $0<y<1$, it contains a term $\ln y$, and the position of the front a term $\ln(1/\delta)$, exactly when $\ell$ divides $q-1$. For $\ell=1$ this is always the case, while for $y'=y^2-y^4$ there is no logarithmic shift.

\subsection{Phase delays of fixed-step methods}
Theorem~\ref{prop:delay} applies with $q=3$; the densities, obtained symbolically, are
\begin{gather*}
g_{\mathrm{Euler}}=\frac{3-4y}{2y(1-y)}=-g_{\mathrm{impl.\,Euler}},\qquad g_{\mathrm{Heun}}=\frac{y(10y^2-15y+6)}{6(1-y)},\\
g_{\mathrm{RK4}}=\frac{y^5(15-30y-3y^2+30y^3-10y^4)}{240(1-y)} .
\end{gather*}
Hence $\kappa=\pm\frac32$ for the Euler methods ($=q/2$), while for Heun and RK4 the density vanishes at $y=0$ (like $y^{p(q-1)-q}$), even faster than for $q=2$. Fig.~\ref{fig:q3}(a) and the lower half of Table~\ref{tab:remainder} confirm this. The Euler delays grow by $\frac32h\ln10=0.86$ per decade of $\delta$ for $h=0.25$. The Heun delay is $(4.1$--$4.2)\cdot10^{-3}$ and the RK4 delay $(2.7$--$3.0)\cdot10^{-7}$ for every $\delta\in[10^{-4},10^{-2}]$, although the number of steps grows from $2\cdot10^4$ to $2\cdot10^8$. The resolved scaled remainders are again consistent with a bound uniform in $\delta$.

\subsection{Adaptive solvers}
\begin{figure}[t]
\centering
\includegraphics[width=\textwidth]{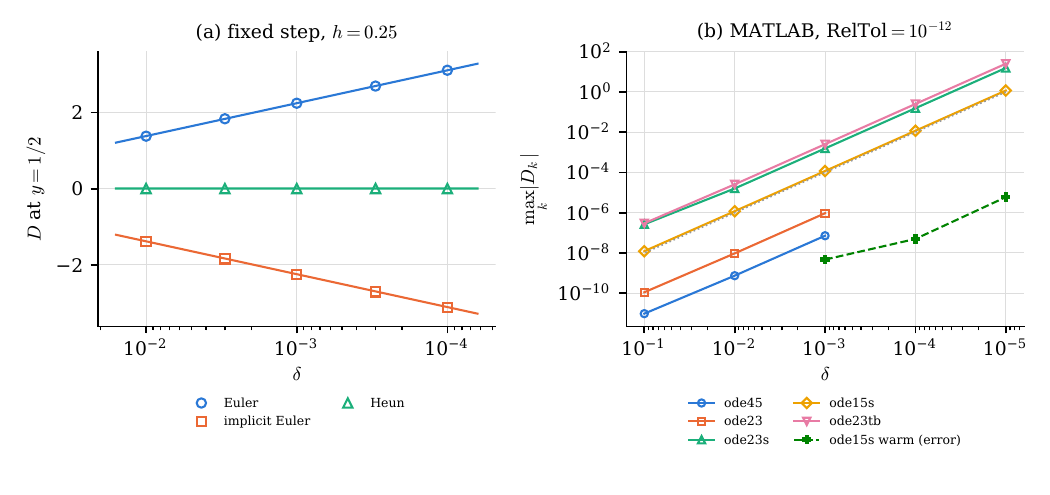}
\caption{The case $q=3$. (a) Phase delays at $y=1/2$, $h=0.25$ (markers: exact; lines: Theorem~\ref{prop:delay}). (b) Maximal phase delay of the MATLAB solvers at $\mathrm{RelTol}=10^{-12}$, $\mathrm{AbsTol}=10^{-20}$, and error of warm-started \texttt{ode15s}; the dotted line has slope $\delta^{-2}$.}
\label{fig:q3}
\end{figure}

The argument of Section~\ref{sec:adaptive} now gives phase increments $\mathrm{RelTol}\cdot y/f_3(y)\approx\mathrm{RelTol}/y^2$ per step, hence $D\approx K\,\mathrm{RelTol}^{\alpha}/\delta^{2}$. The MATLAB runs on $t\in[0,1/\delta^2]$ (Fig.~\ref{fig:q3}b, Table~S5) show exactly this. In the relative regime, $\mathrm{RelTol}=10^{-12}$, the delay of \texttt{ode15s} grows from $1.2\cdot10^{-8}$ at $\delta=10^{-1}$ to $1.2$ at $\delta=10^{-5}$, a factor of $100$ per decade (the other solvers: $59$--$133$), and at the benchmark tolerance the stiff solvers have errors above $0.94$ already for $\delta=10^{-3}$. The median ratio of maximal error to maximal delay over the completed runs with $\max|D_k|<1$ is $0.105$, equal to $m_3=27/256$, as Lemma~\ref{lem:phase} predicts. The warm start from $y=0.05$ again removes the $\delta$-dependence: at $\mathrm{RelTol}=10^{-8}$, \texttt{ode15s} reaches an error of $1.1\cdot10^{-5}$ with $332$--$348$ steps for every $\delta\in[10^{-5},10^{-3}]$, against errors above $0.99$ at $\delta=10^{-5}$ without it.

%=====================================================================
\section{Conclusions}

Every solution of the flame problem is a translate of one profile. Made quantitative, this observation yields the following.

\paragraph{Proved (for $y'=y^q(1-y)$, $q\ge2$)}
\begin{itemize}
\item The global error of every numerical solution with values in $[\delta,1)$ is exactly a phase error (Lemma~\ref{lem:phase}), and for an arbitrary step sequence the delay is the sum of the local errors converted into time (Proposition~\ref{prop:defects}).
\item For fixed-step one-step methods of order $p$ the delay is $h^p\int_\delta^{y}g$ with a remainder $O(h^{p+1})$ uniform in $\delta$ (Theorem~\ref{prop:delay}): it grows at most like $h\ln(1/\delta)$ for $p=1$ (exactly so for the Euler methods) and stays bounded for $p\ge2$. The proof uses only the scaling structure of the method near $y=0$.
\end{itemize}

\paragraph{Modeled and observed (MATLAB R2023b; SciPy, SUNDIALS and Julia for the remedies)}
\begin{itemize}
\item The step counts are predictable: an explicit Runge--Kutta code spends about $L_{\mathrm{stiff}}/\beta$ steps in the stiff phase (within $2\%$ for \texttt{ode45} and \texttt{ode23}; $3022$ predicted, $3007$ observed for \texttt{ode45} at $\delta=10^{-4}$), while the slow phase costs every adaptive code only $O(\ln(1/\delta))$ steps and is where the stiff codes spend almost all of theirs.
\item The adaptive codes incur a delay of the order $\mathrm{RelTol}^{\alpha}/\delta^{q-1}$ (or $\mathrm{AbsTol}^{\alpha}/\delta^{q-1+\alpha}$), $\alpha=p/(p+1)$ or $1$, committed in the first few steps. This law is a heuristic model; the measured exponents agree with it. At the benchmark tolerance the fronts computed by \texttt{ode23s} and \texttt{ode15s} for $\delta=10^{-4}$ are $29$ and $329$ time units early, with errors close to $1$.
\item Tightening the tolerances reduces the delay, but the tolerance must fall by about two decades or more per decade of $\delta$, and for the stiff codes the cost and the limit $r_{\min}$ make this impractical below $\delta\approx10^{-4}$.
\item Controlling time with an absolute tolerance reduces the error strongly. The warm start from the exact translate removes the $\delta$-dependence; time as a dependent variable reduces it, with a residual dependence on $\delta$ for the low-order codes; and integrating for the clock $w=y^{1-q}$ with an absolute tolerance kept the maximal nodal error below $4\cdot10^{-3}$ at $\mathrm{tol}=10^{-4}$ and $6\cdot10^{-5}$ at $\mathrm{tol}=10^{-8}$ for every $\delta$ tested, at a median cost close to that of the standard formulation, within the limits set by $r_{\min}$ and by rounding (on a common output grid the largest error found is $2.5\cdot10^{-2}$, for \texttt{ode23t}).
\end{itemize}

\paragraph{Beyond the flame family} Nothing in the mechanism is special to the flame problem: in every autonomous problem a local error moves the solution along its orbit, and the flow does not correct the shift. The general theory, with a condition number of the event time for systems, a detection algorithm and applications to an epidemic model, hydrogen--air ignition, neuron models and a heteroclinic cycle, is given in \cite{BorosC2026}. For the flame family itself two questions remain open: a proof that the remainder of Theorem~\ref{prop:delay} is uniform in $\delta$ also for the transformed equations of Remark~\ref{rem:variables}, and a global error bound for the clock formulation.

\section*{Acknowledgements}
This work grew out of a problem proposed by Emil C\u{a}tina\cb{s} about ten years ago at the Tiberiu Popoviciu Institute of Numerical Analysis. The author thanks him for introducing the problem to him, for his constant encouragement, and for many stimulating discussions, and thanks Mihai Nechita for the earlier collaboration on this topic.

\section*{CRediT authorship contribution statement}
\textbf{Imre Boros:} Conceptualization, Methodology, Formal analysis, Software, Investigation, Validation, Visualization, Writing -- original draft, Writing -- review \& editing.

\section*{Declaration of competing interest}
The author declares that he has no known competing financial interests or personal relationships that could have appeared to influence the work reported in this paper.

\section*{Funding}
This research did not receive any specific grant from funding agencies in the public, commercial, or not-for-profit sectors.

\section*{Data availability}
The MATLAB, Python, C and Julia code and the data that reproduce all figures and tables, and the supplementary Tables~S1--S6, are openly available at \url{https://github.com/boros91/flame-phase-errors}; the version used in this paper (v1.0.0) is archived on Zenodo, \url{https://doi.org/10.5281/zenodo.24314322} \cite{BorosCode2026}.

\section*{Declaration of generative AI and AI-assisted technologies in the manuscript preparation process}
During the preparation of this work the author used the AI assistant Claude (Anthropic) in order to help write and check the code, verify the numerical results against the saved outputs, and edit the text. After using this tool, the author reviewed and edited the content as needed and takes full responsibility for the content of the published article.

%=====================================================================


\begin{thebibliography}{99}\sloppy

\bibitem{BergerKohn1988}
M. Berger, R.V. Kohn, A rescaling algorithm for the numerical calculation of blowing-up solutions, Commun. Pure Appl. Math. 41 (1988) 841--863.

\bibitem{BorosB2026}
I. Boros, Front-pinned spectral collocation for the flame problem and related ignition problems, preprint, arXiv:2610.27864 (2026), \url{https://doi.org/10.48550/arXiv.2610.27864}.

\bibitem{BorosC2026}
I. Boros, Phase ill-conditioning of autonomous initial value problems: event times under local error control, manuscript, 2026.

\bibitem{BorosCode2026}
I. Boros, flame-phase-errors: code and data for ``Right shape, wrong time: phase errors of ODE solvers in the flame problem'', version 1.0.0, Zenodo, 2026, \url{https://doi.org/10.5281/zenodo.24314322}.

\bibitem{BuddPiggott2001}
C.J. Budd, M.D. Piggott, The geometric integration of scale-invariant ordinary and partial differential equations, J. Comput. Appl. Math. 128 (2001) 399--422.

\bibitem{CanoSanzSerna1997}
B. Cano, J.M. Sanz-Serna, Error growth in the numerical integration of periodic orbits, with application to Hamiltonian and reversible systems, SIAM J. Numer. Anal. 34 (1997) 1391--1417.

\bibitem{ChaudhryEstepStevensTavener2021}
J.H. Chaudhry, D. Estep, Z. Stevens, S.J. Tavener, Error estimation and uncertainty quantification for first time to a threshold value, BIT Numer. Math. 61 (2021) 275--307, \url{https://doi.org/10.1007/s10543-020-00825-0}.

\bibitem{Corless1996}
R.M. Corless, G.H. Gonnet, D.E.G. Hare, D.J. Jeffrey, D.E. Knuth, On the Lambert $W$ function, Adv. Comput. Math. 5 (1996) 329--359.

\bibitem{CorlessJeffrey2002}
R.M. Corless, D.J. Jeffrey, The Wright $\omega$ function, in: Artificial Intelligence, Automated Reasoning, and Symbolic Computation, Lecture Notes in Comput. Sci., vol. 2385, Springer, Berlin, 2002, pp. 76--89.

\bibitem{GrayScott1990}
P. Gray, S.K. Scott, Chemical Oscillations and Instabilities: Non-linear Chemical Kinetics, Oxford University Press, Oxford, 1990.

\bibitem{Hairer1997}
E. Hairer, Variable time step integration with symplectic methods, Appl. Numer. Math. 25 (1997) 219--227.

\bibitem{HairerNorsettWanner1993}
E. Hairer, S.P. N{\o}rsett, G. Wanner, Solving Ordinary Differential Equations I. Nonstiff Problems, second ed., Springer, Berlin, 1993.

\bibitem{HairerWanner1996}
E. Hairer, G. Wanner, Solving Ordinary Differential Equations II. Stiff and Differential-Algebraic Problems, second ed., Springer, Berlin, 1996.

\bibitem{HallHigham1988}
G. Hall, D.J. Higham, Analysis of stepsize selection schemes for Runge--Kutta codes, IMA J. Numer. Anal. 8 (1988) 305--310.

\bibitem{Higham1993}
D.J. Higham, The tolerance proportionality of adaptive ODE solvers, J. Comput. Appl. Math. 45 (1993) 227--236, \url{https://doi.org/10.1016/0377-0427(93)90277-I}.

\bibitem{HighamHumphriesWain2000}
D.J. Higham, A.R. Humphries, R.J. Wain, Phase space error control for dynamical systems, SIAM J. Sci. Comput. 21 (2000) 2275--2294.

\bibitem{Hindmarsh2005}
A.C. Hindmarsh, P.N. Brown, K.E. Grant, S.L. Lee, R. Serban, D.E. Shumaker, C.S. Woodward, SUNDIALS: suite of nonlinear and differential/algebraic equation solvers, ACM Trans. Math. Softw. 31 (2005) 363--396.

\bibitem{HoelMeo2026}
H. Hoel, J.V. Meo, A priori adaptive numerical methods for estimating blow-up times of autonomous ODEs, preprint, arXiv:2603.12957 (2026), \url{https://doi.org/10.48550/arXiv.2603.12957}.

\bibitem{HumphriesVigneswaran2004}
A.R. Humphries, R. Vigneswaran, Phase space stability error control with variable time-stepping Runge--Kutta methods for dynamical systems, Appl. Numer. Anal. Comput. Math. 1 (2004) 469--488, \url{https://doi.org/10.1002/anac.200410010}.

\bibitem{LuLawYooChen2009}
T. Lu, C.K. Law, C.S. Yoo, J.H. Chen, Dynamic stiffness removal for direct numerical simulations, Combust. Flame 156 (2009) 1542--1551.

\bibitem{MeisrimelBirken2020}
P. Meisrimel, P. Birken, Goal oriented time adaptivity using local error estimates, Algorithms 13 (2020) 113, \url{https://doi.org/10.3390/a13050113}.

\bibitem{Moler2004}
C. Moler, Numerical Computing with MATLAB, SIAM, Philadelphia, 2004.

\bibitem{Moler2014}
C. Moler, Ordinary differential equations, stiffness, Cleve's Corner (MathWorks blog), June 9, 2014, \url{https://blogs.mathworks.com/cleve/2014/06/09/ordinary-differential-equations-stiffness/} (accessed 26 October 2026).

\bibitem{MottOranvanLeer2000}
D.R. Mott, E.S. Oran, B. van Leer, A quasi-steady-state solver for the stiff ordinary differential equations of reaction kinetics, J. Comput. Phys. 164 (2000) 407--428.

\bibitem{Nechita2016}
M. Nechita, Revisiting a flame problem. Remarks on some non-standard finite difference schemes, Didact. Math. 34 (2016) 51--56.

\bibitem{OzkayaGaugerNemili2016}
E. \"Ozkaya, N.R. Gauger, A. Nemili, Chaotic behavior of stiff ODEs and their derivatives: an illustrative example, preprint, arXiv:1610.03358 (2016), \url{https://doi.org/10.48550/arXiv.1610.03358}.

\bibitem{Petzold1983}
L.R. Petzold, Automatic selection of methods for solving stiff and nonstiff systems of ordinary differential equations, SIAM J. Sci. Stat. Comput. 4 (1983) 136--148.

\bibitem{RackauckasNie2017}
C. Rackauckas, Q. Nie, DifferentialEquations.jl -- a performant and feature-rich ecosystem for solving differential equations in Julia, J. Open Res. Softw. 5 (2017) 15.

\bibitem{Shampine1977}
L.F. Shampine, Stiffness and nonstiff differential equation solvers, II: Detecting stiffness with Runge--Kutta methods, ACM Trans. Math. Softw. 3 (1977) 44--53.

\bibitem{ShampineGladwellThompson2003}
L.F. Shampine, I. Gladwell, S. Thompson, Solving ODEs with MATLAB, Cambridge University Press, Cambridge, 2003.

\bibitem{ShampineReichelt1997}
L.F. Shampine, M.W. Reichelt, The MATLAB ODE suite, SIAM J. Sci. Comput. 18 (1997) 1--22.

\bibitem{ShampineThompson2007}
L.F. Shampine, S. Thompson, Stiff systems, Scholarpedia 2 (2007) 2855.

\bibitem{SoderlindWang2006}
G. S{\"o}derlind, L. Wang, Evaluating numerical ODE/DAE methods, algorithms and software, J. Comput. Appl. Math. 185 (2006) 244--260, \url{https://doi.org/10.1016/j.cam.2005.03.009}.

\bibitem{StuartFloater1990}
A.M. Stuart, M.S. Floater, On the computation of blow-up, Eur. J. Appl. Math. 1 (1990) 47--71.

\bibitem{Sundman1912}
K.F. Sundman, M\'emoire sur le probl\`eme des trois corps, Acta Math. 36 (1912) 105--179.

\bibitem{SciPy2020}
P. Virtanen, R. Gommers, T.E. Oliphant, et al., SciPy 1.0: fundamental algorithms for scientific computing in Python, Nat. Methods 17 (2020) 261--272.

\end{thebibliography}
\end{document}